\documentclass[12pt]{article}
\usepackage[all]{xy}
\usepackage{amsmath}
\usepackage{amsfonts}
\usepackage{amssymb}
\usepackage{amscd}
\usepackage{amsthm}
\usepackage{latexsym}
\usepackage{amsbsy}
\newtheorem{lem}{Lemma}[section]
\newtheorem{prop}{Proposition}[section]
\newtheorem{theorem}{Theorem}[section]
\newtheorem{corollary}{Corollary}[section]
\newtheorem{definition}{Definition}[section]
\newtheorem{example}{Example}[section]
\newtheorem{Conjecture}{Conjecture}[section]

\newcommand{\de}{\delta}
\newcommand{\ot}{\otimes}
\begin{document}
\title{Cyclic  Cohomology of  Hopf  algebras and Hopf Algebroids}
\author {M. Khalkhali,~~~ B. Rangipour,
\\\texttt{~masoud@uwo.ca ~~~~brangipo@uwo.ca}
\\Department of Mathematics 
\\University of Western Ontario }
\maketitle

\begin{abstract}
We review recent progress in the study of cyclic cohomology of Hopf algebras, 
 Hopf algebroids, and  invariant cyclic homology starting with the
pioneering work of Connes-Moscovici.
 \end{abstract}
 \section{Introduction} 
 It is well known that the theory of characteristic classes of vector bundles, 
 more precisely the Chern character, can
  be extended to the noncommutative geometry,
  thanks to the noncommutative Chern-Weil theory of Connes \cite{ac80,ac,ac88}.
   In order to have a similar extension for
   quantum principal bundles, for example Hopf-Galois extensions, 
    one needs first appropriate analogues of group 
   and Lie algebra cohomology  of  Hopf algebras. The recent works of Connes-Moscovici ~\cite{achm98,achm00,achm01}
   on the index theory of  transversely elliptic operators,
   more precisely their  definition of cyclic cohomology of Hopf algebras, 
   provides one with such a theory.
   
   It is the goal of the present article to review the developments in the study
    of cyclic cohomology of Hopf algebras, starting  with the  pioneering work of 
	Connes-Moscovici ~\cite{achm98,achm00,achm01}. We will present a dual cyclic 
	theory for Hopf algebras, first  defined in \cite{kr}, and independently in 
	 ~\cite{ta}. One motivation is that, as it was observed by M. Crainic
    ~\cite{cr},
    cyclic cohomology of  cosemisimple Hopf algebras, { e.g.} the algebra of polynomial functions on a compact  
     quantum 	groups, due to   existence of Haar integral, is always
     trivial. In other words it behaves in much the same way as continuous group cohomology. 
    Let $HP^\bullet$ and $\widetilde {HP}_\bullet$ denote the resulting periodic 
	cyclic (co)homology  groups in the sense of ~\cite{achm98} and ~\cite{kr}, respectively. We present
     two 
    very general results: for any commutative Hopf algebra $\mathcal{H}$,   
    	 $HP^\bullet({\mathcal{H}})$ decomposes into direct  sums of Hochschild 
	 cohomology groups of the coalgebra $\mathcal{H}$ with trivial coefficients, 
	 and for any cocommutative $\mathcal{H}$,   $\widetilde {HP}_\bullet(\mathcal{H})$
	  decomposes as Hochschild homology groups of algebra $\mathcal{H}$ with trivial 
	  coefficients. So far very few examples of computations of $HP^\bullet$ and
      $\widetilde {HP}_\bullet$  for quantum groups are known. We present what 
	  is known  in Sections 3 and 4.
	        
       In Section  5 we review the main results on cyclic cohomology of extended 
	   Hopf algebras known so far, following \cite{achm01,kr2}. Extended Hopf 
	   algebras are closely related to Hopf algebroids. It seems that  now the 
       question of finding an appropriate algebraic framework to define cyclic 
	   cohomology of Hopf algebroids is settled by ~\cite{kr2}.
       
       In Section 6 we present some of the results obtained in \cite{ak} on cyclic
        cohomology of smash products. 
        
Cyclic homology of Hopf algebras can be understood from two distinct     points of view. The first view, 
 due to Connes and Moscovici  \cite{achm00,achm99,achm98}, is based on the existence of characteristic map for (co)action
 of Hopf algebras on algebras (see the introductory remarks in Section \ref{sec4} 
 for more on this). In the second point of view, first advocated in \cite{kr5}, cyclic (co)homology 
 of Hopf algebras appears as a special case  of a more general theory called invariant cyclic homology. We review 
 this theory in Section \ref{sec7}.  It turns out that the invariant cyclic homology of Hopf algebra is isomorphic   
  to its Hopf algebraic cyclic homology. This is  remarkably similar to interpreting the cohomology 
 of the Lie algebra of a 
 Lie group as invariant de Rham cohomology of its Lie group as is done  by Chevalley and Eilenberg \cite{ce}.
            
        It was not our intention to cover all aspects of this new branch of noncommutative geometry in this paper.
         For applications to transverse index theory and for the whole theory one should consult the original 
         Connes-Moscovici articles 
          ~\cite{achm98,achm99,achm01}   as well as their  review article~\cite{achm00}. We also 
           recommend ~\cite{va} for a  general introduction to applications of Hopf algebras 
           in noncommutative geometry. Much remains to be done in this area. For example, the relation between 
            cyclic homology of Hopf algebras and developments in Hopf-Galois theory (see { e.g. Montgomery's 
            book~\cite{mo}}) remain to be explored.
             Also, what is missing is a general conjecture about the nature of Hopf cyclic homology of the algebra 
             of polynomial functions (or smooth functions, provided they are defined) of  quantum groups. 
      
        \section{Preliminaries on Hopf algebras }
In this paper algebra means an associative, not necessarily commutative, 
unital algebra over a fixed commutative ground ring $k$.
Similar convention applies to \textit{coalgebras}, \textit{bialgebras} and \textit{Hopf algebrs}.
 The undecorated tensor product
$\otimes$ means the tensor product over $k$.   If $\mathcal{H}$ is a 
Hopf algebra, we denote its coproduct by $\Delta: \mathcal{H}\longrightarrow \mathcal{H}\otimes \mathcal{H}$, 
its counit by $\epsilon : \mathcal{H}\longrightarrow k$, its unit by $\eta : k\longrightarrow \mathcal{H}$ 
and its antipode by $S: \mathcal{H}\longrightarrow \mathcal{H}$. We will use Sweedler's
 notation $\Delta(h)=h^{(1)}\otimes h^{(2)}$, $(\Delta\otimes id)\Delta(h)=h^{(1)}\otimes h^{(2)}\otimes h^{(3)}$,
 etc, where summation is understood.

If $\mathcal{H}$ is a Hopf algebra, the word $\mathcal{H}$-module means a module over the 
underlying algebra of $\mathcal{H}$. 
Similarly, an $\mathcal{H}$-comodule is a comodule over the underlying coalgebra of $\mathcal{H}$. 
The same convention applies to $\mathcal{H}$-bimodules and $\mathcal{H}$-bicomodules. 
 The category of (left) $\mathcal{H}$-modules has a tensor product defined via the coprouct of $\mathcal{H}$:
if $M$ and $N$ are left $\mathcal{H}$-modules, their
 tensor product $M\otimes N$ is again an $\mathcal{H}$-module via $$h(m\otimes n)= h^{(1)}m\otimes h^{(2)}n.$$
similarly, if $M$ and $N$ are left $\mathcal{H}$-comodules, the tensor product $M\otimes N$ 
is again an $\mathcal{H}$-comodule 
via $$\Delta (m\otimes n)=m^{(-1)} n^{(-1)}\otimes m^{(0)}\otimes n ^{(0)}.$$

We take the point of view, standard in noncommutative geometry, that a noncommutative space is encoded by 
an algebra or by a coalgebra. The idea of \textit{symmetry,  i.e}. action of a group on a space, 
can be expressed by the action/coaction of a Hopf algebra on an algebra/coalgebra. Thus four possibilities arise.
 Let $\mathcal{H}$ be a Hopf algebra. An algebra $A$ is called a left $\mathcal{H}$-\textit{module algebra} 
 if it  is a left $\mathcal{H}$-module and the multiplication  map   $A\otimes A\longrightarrow A$ and the unit map
 are   morphisms of  $\mathcal{H}$-modules. That is $$h(ab)=h^{(1)}(a)h^{(2)}(b),\qquad h(1)=\epsilon(h)1,$$
for $h\in \mathcal{H} ,a,b\in A$. Similarly an algebra $A$ is called a $\mathcal{H}$-comodule algebra, if  
$A$ is a left $\mathcal{H}$-comodule and the multiplication and the unit maps 
 are  morphisms of $\mathcal{H}$-comodules. 
In a similar fashion an  $\mathcal{H}$-module coalgebra is a coalgebra $C$ which is
 a left $\mathcal{H}$-module, and the comultiplication $\Delta: C\longrightarrow C\otimes C$ and the counit map are 
 $\mathcal{H}$-module maps. Finally  an $\mathcal{H}$-comodule 
 coalgebra is a coalgebra $C$ which is an $\mathcal{H}$-comodule 
and the coproduct and counit map are comodule maps.

The \textit{smash product} $A\# \mathcal{H}$ of an $\mathcal{H}$-module
 algebra $A$ with $\mathcal{H}$ is,  as a $k$-module,  $A\otimes \mathcal{H}$  with the product
$$(a\otimes g)(b\otimes h)=a (g^{(1)}b)\otimes g^{(2)}h.$$
It is an associative algebra under the above product.\\
{\bf Examples}
\begin{itemize}
\item{1.} For $\mathcal{H}=U(\mathfrak{g})$, the enveloping algebra of a Lie algebra, $A$ is 
an $\mathcal{H}$-module algebra     
iff  $\mathfrak{g}$ acts on $A$ by derivations, { i.e.}  we have a Lie algebra map
 $\mathfrak{g}\longrightarrow Der(A)$.\\    
\item{2.} For $\mathcal{H}=kG$, the group algebra of a (discrete) group $G$, $A$ is a $\mathcal{H}$-module
 algebra iff  $G$ acts on $A$ via automorphisms $G\longrightarrow Aut(A)$. The smash product $A\# \mathcal{H}$ is
 then isomorphic to the crossed product algebra $A\rtimes G$.\\
\item{3.} For any Hopf algebra $\mathcal{H}$, the algebra $A=\mathcal{H}$ is 
an $\mathcal{H}$-comodule algebra where the coaction  is afforded by comultiplication
 $\mathcal{H}\longrightarrow \mathcal{H}\otimes \mathcal{H}$. Similarly, the coalgebra $\mathcal{H}$ is an
  $\mathcal{H}$-module coalgebra where the action is given by the multiplication 
  $\mathcal{H}\otimes\mathcal{H}\rightarrow \mathcal{H}$. These are analogues  of the action  of a group
   on itself by translations.
\item{4.} By a theorem of  Kostant ~\cite{sw},  any cocommutative Hopf algebra $\mathcal{H}$ over an algebraically 
closed field of characteristic zero is isomorphic (as  a Hopf algebra ) with a 
  smash product  $\mathcal{H}= U(P(\mathcal{H}))\# kG(\mathcal{H})$,  
where  $P(\mathcal{H})$ is the  Lie algebra  of  primitive elements of $\mathcal{H}$ and $G(\mathcal{H})$ is the
 group of all grouplike
 elements of $\mathcal{H}$ and $G(\mathcal{H})$ acts on $P(\mathcal{H})$ by inner 
 automorphisms $(g,h)\mapsto ghg^{-1}$, for $g\in G(\mathcal{H})$ and $h\in P(\mathcal{H})$.
\end{itemize}

\section{Cyclic modules}
Cyclic co/homology was first defined for (associative) algebras through explicit 
complexes or bicomplexes. Soon after,  Connes introduced the notion of cyclic
 module and defined cyclic homology of cyclic modules \cite{ac88}. The motivation was 
 to define cyclic homology of algebras as a derived functor. Since the category 
 of algebras and algebra homomorphisms is not an additive category, the standard 
(abelian) homological algebra is not enough. In Connes' approach,  the category of
 cyclic modules appears as \textquoteleft\textquoteleft abelianization" of the category of
algebras with the embedding defined by the functor $A \mapsto A^\natural$, 
explained  below. For an alternative approach one can consult (\cite{ft}), where 
cyclic cohomology is shown to be the nonabelian derived functor of the functor 
of traces on $A$.   It was soon realized that cyclic modules and the flexibility 
they afford are indispensable tools in the theory. A recent example is the cyclic
 homology of Hopf algebras which can not be defined as the cyclic homology of 
an algebra or coalgebra.

In this section we recall  the theory of cyclic and paracyclic modules and their 
cyclic  homologies. We also consider the doubly graded version i.e.
 biparacyclic modules and the generalized Eilenberg-Zilber theorem~\cite{ac88,ft,gj}.

For $r \geq 1$ an integer or $r=\infty$, let $\Lambda^r$ denote 
the r-cyclic category.  An {\it r-cyclic object} in a category $\mathcal{C}$ is 
a contravariant functor $\Lambda^r \rightarrow \mathcal{C}$. Equivalently, we have 
a sequence $X_n, n \geq 0$, of objects of $\mathcal{C}$ and morphisms called face, 
 degeneracy and cyclic operators 
 $$\delta_i: X_{n} \rightarrow X_{n-1},\quad  \sigma_i : X_n \rightarrow X_{n+1}, 
 \quad  \tau : X_n \rightarrow X_n \qquad  0\le i\le n$$
 such that $(X,\delta_i,\sigma_i)$ is a simplicial object and the following extra 
  relations are satisfied:
\begin{eqnarray*}
\delta_i \tau&=& \tau\delta_{i-1}\hspace{43 pt} 1\le i\le
n\\ 
\delta_0 \tau&=& \delta_{n}
\\ \sigma_i \tau&=& \tau\sigma _{i-1} \hspace{43 pt} 1\le i\le n\\
\sigma_0 \tau&=& \tau^2\sigma_n \\   
\tau^{r(n+1)} &=& \mbox{id} _n . 
\end{eqnarray*}
For $r=\infty$, the last relation is replaced by the empty relation and we have 
 a  {\it paracyclic} object. For $r=1$, a   $ \Lambda^1$   object  is  a  
 {\it cyclic object}.

A cocyclic object is defined in a dual manner. Thus a cocyclic object in 
 $\mathcal{C}$ is a covariant functor $\Lambda^1 \rightarrow \mathcal{C}$. Let $k$
  be a commutative ground ring.  A cyclic module over $k$ is a cyclic object in 
  the category of $k$-modules. We denote the category of cyclic $k$-modules by $\Lambda_k $.
  
   Next,  let us recall that a {\it biparacyclic } 
   object  in a category $\mathcal{C}$ is a contravariant functor $\Lambda^\infty
  \times\Lambda^\infty \rightarrow \mathcal{C}$. Equivalently, we have a doubly 
  graded set of objects $X_{n,m}$, $n,m\geq 0$ in $\mathcal{C}$ with  horizontal and vertical 
    face,  degeneracy and cyclic operators  $\delta_i,\sigma_i,\tau, d_i, s_i,t  $
   such  that each row and  each column  is a paracyclic object in $\mathcal{C}$
    and vertical and horizontal operators  commute.  A biparacyclic object $X$ is called {\it cylindrical}
	 if the operators $\tau^{m+1},t^{n+1}:X_{m,n}\rightarrow X_{m,n}$  are
    inverse of each other.
If $X$ is cylindrical then it is easy to see that its {\it diagonal},  $d(X)$, 
defined by $d(X)_n=X_{n,n} $ with face,  degeneracy and cyclic maps $\delta _i d_i$,
 $\sigma _is_i$ and  $\tau  t$ is a cyclic object. 
  
We give a few examples of cyclic modules that will be used in this paper. The
 first example is the most fundamental example which motivated the whole theory.
\begin{itemize}
\item[1.]  Let $A$ be an algebra. The cyclic module $A^ \natural$ 
is defined by $A_n^\natural=A^{\otimes(n+1)} , n\geq 0$, with the face,  degeneracy 
and cyclic operators  defined by
\begin{eqnarray*}
\delta_i(a_0 \otimes a_1\otimes \dots \otimes a_n)&=&a_0 \otimes \dots \otimes 
a_{i}a_{i+1}\otimes \dots \otimes a_n\\
\delta_n(a_0 \otimes a_1\otimes\dots \otimes a_n)&=&a_na_0 \otimes a_1 \otimes  
\dots  \otimes a_{n-1}\\
\sigma_i(a_0 \otimes a_1\otimes\dots \otimes a_n)&=&a_0 \otimes\dots \otimes  
a_i \otimes  1 \otimes \dots\otimes
 a_n\\
\tau(a_0 \otimes a_1\otimes\dots \otimes a_n)&=& a_n \otimes a_0\dots \otimes a_{n-1}.
\end{eqnarray*}
The underlying simplicial module of $A^\natural $ is a special case of the following
 simplicial module. Let $M$ be an $A$-bimodule. Let $C_n(A,M)=M \otimes A^{\otimes n},\quad 
 n\geq 0$.  For  $n=0$, 
we put $C_0(A,M)=M$. Then the following faces and degeneracies $\delta _i, 
\sigma _i$ define a simplicial module structure on $C_\bullet(A,M)$:
\begin{eqnarray*}
\delta _0(m \otimes a_1 \otimes  \dots \otimes a_n)&=& m a_1 \otimes a_2 \otimes
\dots \otimes a_n\\
\delta _i(m \otimes a_1 \otimes \dots \otimes a_n)&=& m \otimes a_1 \otimes 
\dots \otimes a_{i}a_{i+1}\otimes \dots
\otimes a_n\\
\delta_n(m \otimes a_1 \otimes \dots \otimes a_n)&=& a_n m \otimes a_1 \otimes
 \dots \otimes a_{n-1}\\
 \sigma _0(m \otimes a_1 \otimes \dots  \otimes a_n)&=& m\otimes 1\otimes a_1 \otimes   
\dots \otimes a_n\\
\sigma _i(m \otimes a_1 \otimes \dots  \otimes a_n)&=& m\otimes a_1  \otimes  
\dots \otimes  a_i \otimes 1 \otimes\dots \otimes a_n \quad 1\le i\le n.
\end{eqnarray*}
Obviously, for $M=A$ we obtain $A^\natural$. In general, there is no cyclic 
structure on $C_\bullet{(A,M)}$.\\
\item[2.]  let $C$ be a coalgebra. The cocyclic module $C_\natural$ is 
defined by $C_\natural ^n=C^{\otimes n+1},\;\; n\geq 0$, with coface, codegeneracy
 and cyclic operators:
\begin{eqnarray*}
\delta_i(c_0 \otimes c_1 \otimes \dots \otimes c_n)&=& c_0 \otimes\dots 
\otimes  c_i^{(1)}\otimes c_i^{(2)}\otimes c_n ~~~ 0 \leq i \leq n\\
\delta_{n+1}(c_0 \otimes c_1 \otimes \dots \otimes c_n)&=& c_0^{(2)}\otimes c_1 
\otimes \dots \otimes c_n  \otimes c_0^{(1)}\\
\sigma_i(c_0 \otimes c_1 \otimes \dots \otimes c_n)&=&c_0 \otimes \dots 
c_i \otimes \varepsilon(c_{i+1})\otimes 
\dots\otimes c_n ~~~0\leq i \leq n-1\\
\tau(c_0 \otimes c_1 \otimes \dots \otimes c_n)  &=& c_1 \otimes c_2 \otimes 
\dots \otimes c_n \otimes c_0,
\end{eqnarray*}
where as usual $\Delta(c)= c^{(1)}\otimes c^{(2)} $ (Sweedler's notation).
The underlying cosimplicial module for $C_\natural$ is a special case of the 
following cosimplicial module. Let $M$ be a $C$-bicomodule and $C^n(C,M)= M 
\otimes C^{\otimes n}$.
The following coface and codegeneracy operators define a cosimplicial module.
\begin{eqnarray*}
\delta_0(m \otimes c_1 \otimes \dots \otimes c_n)&=&  m^{(0)}\otimes m^{(1)}\otimes
 c_1 \dots \otimes c_n \\
\delta_i(m \otimes c_1 \otimes \dots \otimes c_n)&=& m \otimes c_1 \dots \otimes  
c_i^{(0)}\otimes c_i^{(1)}\otimes c_n \;\;\text{for }\;\; 1 \leq i \leq n\\
\delta_{n+1}(m\otimes c_1 \otimes \dots \otimes c_n)&=& m_{(0)}\otimes c_1
 \otimes \dots \otimes c_n 
 \otimes m_{(-1)}\\
\sigma_i(m \otimes c_1 \otimes \dots \otimes c_n)&=& m \otimes c_1 \dots 
\varepsilon(c_{i+1})c_{i}\otimes 
\dots\otimes c_n 
~~~~0\leq i \leq n-1, 
\end{eqnarray*}
where we have denoted the left and right comodule maps by $\Delta_l(m)=m_{(-1)}
\otimes m_{(0)}  $ and 
 $\Delta_r(m)=m^{(0)}\otimes m^{(1)}$.   Let $$
 d=\sum_{i=0}^{n+1}(-1)^i\delta_i:C^n(C,M)\rightarrow C^{n+1}(C,M).$$ Then $d^2=0$.
  The cohomology of the 
 complex $(C^\bullet(C,M),d)$  
 is the Hochschild cohomology of the coalgebra $C$ with coefficients in the bicomodule  
 $M$. For $M=C$, we obtain the
  Hochschild  complex of $C_\natural$. Another special case  occurs with $M=k$ and
   $\Delta_r:k \rightarrow k \otimes C \cong C $ and $\Delta _l:k \rightarrow 
   C \otimes k \cong C $, are given by
    $\Delta_r(1)=1\otimes g$ and $\Delta_l(1)=h \otimes 1$, 
   where $g,h\in C$ are grouplike elements. The differential $d :C^{n}\rightarrow 
   C^{n+1}$  in the latter case is given by 
   \begin{multline*}
   d(c_1 \otimes c_2 \dots \otimes c_n)=g \otimes c_1 \otimes  \dots \otimes c_n\\
   +\sum _{i=1}^n(-1)^i c_1 \otimes \dots \otimes \Delta(c_i)\otimes\dots \otimes c_n
   +(-1)^{n+1}c_1\otimes \dots\otimes c_n \otimes h.
   \end{multline*} 
\item[3.] 
 Let $g: A \rightarrow A$ be an  automorphism  of an algebra $A$. 
The paracyclic module $A^\natural _g$ is defined by $A_{g,n}^\natural=A^{\otimes(n+1)}$ 
with the same cyclic structure as $A^\natural$,  except the following changes 
\begin{eqnarray*}
\delta _n(a_0 \otimes a_1 \dots \otimes a_n)&=& g(a_n)a_0 \otimes  \dots \otimes 
a_{n-1}\\
\tau (a_0 \otimes a_1 \dots \otimes a_n)&=& g(a_n)\otimes a_0 \otimes \dots 
 \otimes a_{n-1}.
\end{eqnarray*}
\end{itemize}
One can check that   $A_g^\natural$ is a $\Lambda^\infty $-module and if $g^r=id$, then it 
is a $\Lambda^r$-module. For $g=id$,
 we obtain example 1.

Next, let us indicate how one defines the Hochschild, cyclic and periodic cyclic
 homology of 
a cyclic module. This is particularly important since the cyclic homology of Hopf algebras
 is naturally defined as the cyclic homology of some cyclic  modules  associated with  them.
 Given a cyclic module $M\in \Lambda_k$, its cyclic homology group $HC_n(M)$, $n\ge 0$,
  is defined in (\cite{ac88}) by $$HC_n(M):=Tor_n^{\Lambda_k} (M,k^\natural),$$
  and similarly  the cyclic cohomology groups of $M$ are defined by 
  $$HC^n(M):=Ext^n_{\Lambda_k}(M,k^\natural).$$ 
  
  Using a specific projective resolution  for $k^\natural$, one obtains the following 
  bicomplex to compute  cyclic homology. 
Given a cyclic module $M$, consider the following first quadrant  bicomplex, called
 the {\it cyclic bicomplex}  of  $M$ 
$$\begin{CD}
\vdots @.\vdots @.\vdots @.\\
M_2@<1-\tau<< M_2 @<N<< M_2@<1-\tau<< \dots  \\
@VV bV @VV-b'V @VV bV  \\
M_1@<1+\tau<< M_1 @<N<< M_1@<1+\tau<< \dots \\
@VV bV @VV-b'V @VV bV  \\
M_0@<1-\tau<< M_0 @<N<< M_0@<1-\tau<< \dots 
\end{CD}
$$
We denote this bicomplex by  $CC^+(M)$. The operators $b$, $b'$ and $N$ are  defined by
\begin{eqnarray*}
 b &=& \sum _{i=0}^n(-1)^i \delta _i\\
 b' &=& \sum _{i=0}^{n-1} (-1)^i \delta _i\\
N&=& \sum_{i=0}^n(-1)^{ni} \tau ^i.
\end{eqnarray*}
Using the simplicial and cyclic relations, one can check that $b^2=b'^2=0$, 
$b(1-(-1)^n\tau)=(1-(-1)^{n-1}\tau)b'$
 and $b'N=Nb'$.
 The {\it Hochschild homology} of $M$, denoted $HH_\bullet (M)$, is the homology 
 of the first column $(M_\bullet,b)$.
The {\it  {cyclic homology}} of $M$, denoted by  $H C_\bullet (M)$ is the homology of the total complex
 $Tot C C^+(M)$. 
 
 To define the {\it periodic cyclic} homology of  $M$, 
 we extend the first quadrant  bicomplex $C C^+(M)$ to the left  and denote  it by $C C(M)$. 
  Let $Tot CC(M)$ denote the \textquoteleft\textquoteleft total complex" where instead of direct sums we  use direct product,
   $$TotCC(M)_n=\prod_{i=0}^\infty M_i. $$ It is obviously a 2-periodic complex and its homology is called the
    periodic
    cyclic homology  of $M$ and denoted by  $HP_\bullet(M)$. 
    
The complex $(M_\bullet,b')$ is acyclic with contracting homotopy $\sigma_{-1}=\tau \sigma _n $. One can then show
 that $CC^+(M)$ is homotopy equivalent  to Connes's $(b,B)$ bicomplex 
$$\begin{CD}
\vdots @.\vdots @.\vdots \\
M_2@<B<< M_1@<B<< M_0 \\
@VVbV @VVbV\\
M_1@<B<<M_0\\
@VVbV\\
M_0
\end{CD}
$$ 
Where $B : M_n \rightarrow M_{n+1}$ is Connes's boundary operator defined  by $B=(1-(-1)^n\tau)\sigma_{-1}N$.

 Finally we arrive at the  3rd definition of cyclic homology by noticing that if $k$ 
 is a filed of characteristic zero,
  then the rows of $CC^+(M)$  are acyclic in positive degree and its homology in dimension  zero is
   $$C_n^{\lambda}(M)=\frac{ M_n}{(1-(-1)^n\tau)M_n}.$$
It follows that the total homology, { i.e.}  cyclic homology of $M$ can be computed, if $k$ is a field of 
characteristic zero, as the homology of Connes's cyclic complex $(C_\bullet^\lambda(M),b)$
  
  Now, if $A$ is an associative algebra, its Hochschild, cyclic and periodic cyclic homology, are defined as the 
  corresponding  homology of the cyclic module 
  $A^\natural$.  We denote these groups by $HH_\bullet(A)$, $HC_\bullet(A)$ and $HP_\bullet(A)$, respectively.
  Similarly, if $C$ is a coalgebra, its Hochschild, cyclic and periodic cyclic 
  cohomology are defined as the 
  corresponding homology of the cocyclic module $C_\natural$.        
  
Our next goal is to recall the generalized Eilenberg-Zilber theorem for cylindrical modules from \cite{gj,kr1}.
This is needed in Section $6$ to derive a spectral sequence for cyclic homology of smash products. 

 A {\it parachain complex  } $(M_\bullet,b,B)$ is a chain complex $(M_\bullet,b )$  endowed with a map $B: M_\bullet
  \rightarrow M_{\bullet+1}$
  such that $B^2=0$ and $T=1-(bB+Bb)$ is an invertible operator. For example, a {\it mixed complex } is a 
  parachain  complex such that $bB+Bb=0$. Given a mixed complex $M$ one can define its $(b,B)$-bicomplex  
  as the  Connes' $(b,B)$  bicomplex. One can thus define the Hochschild, cyclic and periodic cyclic homology of mixed
   complexes.
  The definition of {\it bi-parachain complex} should be clear. Given a bi-parachain complex $X_{p,q}$, one defines 
  its total complex $TotX$ by $$(TotX)_n=\oplus X_{p,q},\quad b=b_v+b_h,\quad B=B_v+TB_h,$$ where v and h refers to
   horizontal and vertical differentials. One can check that  $TotX $ is a parachain complex ~\cite{gj}.
  
  Now if $X$ is a cylindrical module and $C(X)$  is the bi-parachain complex obtained by forming the associated 
  mixed complexes 
  horizontally and vertically, then one can check that $Tot(C(X))$ is indeed a mixed complex. On the 
  other hand we know that the diagonal 
  $d(X)$ is a cyclic module and hence its associated chain complex $C(d(X))$ is a mixed complex. 
  
  The following theorem was first proved in ~\cite{gj} using topological arguments. A purely algebraic proof 
  can be found in \cite{kr}.
  
  \begin{theorem}({\cite{gj,kr}})
  Let $X$ be a cylindrical module. There is a quasi-isomorphism of mixed complexes
   $f_0+uf_1:Tot(C(X))\rightarrow C(d(X))$ such that $f_0$ is the shuffle map. 

\end{theorem}

\section{Cyclic cohomology of Hopf algebras}\label{sec4}
Thanks to the  recent work of Connes-Moscovici ~\cite{achm98,achm99,achm01}, the following principle has emerged.
A reasonable co/homology theory for Hopf algebras and Hopf algebra like objects  in noncommutative
 geometry should address the following two issues:
\begin{itemize}
\item[$\bullet$] It should reduce to group co/homology or Lie algebra 
co/homology for $\mathcal{H}=kG$, 
$k\lbrack G\rbrack$ or $U(\mathfrak{g})$
;  Hopf algebras naturally associated to (Lie) groups or Lie algebras.
\item[$\bullet$]  There should exist a characteristic map, connecting the cyclic 
cohomology of a Hopf algebra $\mathcal{H}$ to the  
cyclic cohomology of an algebra
$A$ on which it acts.
For example, for any $\mathcal{H}$-module algebra $A$ and an invariant trace 
$\tau : A \longrightarrow \mathbb{C}$, 
there should exist a  map  
$$\gamma: HC^\bullet(\mathcal{H})\longrightarrow HC^\bullet(A).$$
\end{itemize}

Let us explain both points starting with the first. It might seem 
that given a Hopf algebra $\mathcal{H}$, the 
Hochschild homology of the algebra $\mathcal{H}$
might be a good candidate for  a homology theory for $\mathcal{H}$ in noncommutative geometry. 
After all one knows that for a Lie algebra $\mathfrak{g}$
 and a $U(\mathfrak{g})$-bimodule $M$,
$$H_\bullet(\mathfrak{g}\bold{,}M^{ad})\cong H_\bullet(U(\mathfrak{g}),M)$$ 
where the action of $\mathfrak{g}$ on $M$ is given by $g \cdot m=g m-m g $ \cite{ld}.
Thus Hochschild homology  of $U(\mathfrak{g})$ can be recovered  from the Lie algebra homology 
of $\mathfrak{g}$. Conversely, if $M$ is a $\mathfrak{g}$-module we can turn it into a $U(\mathfrak{g})$-bimodule 
where the left action is induced by $\mathfrak{g}$-action and the right action is by augmentation :  $mX=\epsilon(X)m 
$. It 
follows that  $H_\bullet(\mathfrak{g},M)\cong H_\bullet(U(\mathfrak{g}),M)$, which shows that  the Lie algebra homology 
can also be recovered from Hochschild homology. 
In particular $H_\bullet(\mathfrak{g}\boldsymbol{,}k)\cong H_\bullet(U(\mathfrak{g}),k)$.  Similarly, if $G$ 
is a (discrete)  group and $M$ is a $k G$-bimodule then   
$H_\bullet(G;M^{ad})\cong H H_\bullet(k G,M)$ where the action of  
$G$ on $M^{ad}=M$ {is given by} $ g m=g m g^{-1} $.

In \cite{kr} these type of  results were extended to 
all Hopf algebras in the following way. Let $\mathcal{H}$ be a Hopf algebra and $M$ a left $\mathcal{H}$-module. One 
defines groups ${H}_\bullet(\mathcal{H},M)$ as the left derived functor of the functor of coinvariants
  from $\mathcal{H}$-mod$\rightarrow$$k$-mod,
  $$M\mapsto M_{\mathcal{H}}:=M/\text{ submodule generated by } \lbrace hm - \epsilon(h)m \mid h \in \mathcal{H},\quad
   m \in M \rbrace .$$
   Obviously, $M_{\mathcal{H}}=k \otimes_\mathcal{H} M$ 
  which shows that $H_\bullet(\mathcal{H}\boldsymbol{,} M)\cong Tor_\bullet^\mathcal{H}(k,M)$. For $\mathcal{H}=kG$ or
   $U(\mathfrak{g})$, one obtains group and Lie algebra homologies.   
   
   Now let $\mathcal{H}$ be a Hopf algebra and $M$ be an $\mathcal{H}$-bimodule. We can convert $M$ to a new left 
   $\mathcal{H}$
   -module $M^{ad}=M$, where the action of $\mathcal{H}$ is given by $$h\cdot m=h^{(2)}mS(h^{(1)}).$$ 
 \begin{prop}{(\cite{kr})(Mac Lane isomorphism for Hopf algebras)}\\
Under the above hypotheses there is a canonical isomorphism\\ 
 $$H_n (\mathcal{H} , M) \cong H_n(\mathcal{H} ;
{M^{ad}})= Tor^\mathcal{H}_n(k,M^{ad}),$$  where the left hand side is Hochschild homology.
\end{prop}
Note that the result is true  for all Hopf algebras irrespective of being (co)commutative or not.  

 This suggests to define  $H_\bullet(\mathcal{H}\bold{,}k)$,
where $k$ is an $\mathcal{H}$-bimodule via augmentation map, in analogy with the  group homology. This
is not, however, a reasonable candidate as can be seen by considering $\mathcal{H}=k \lbrack G \rbrack$, 
the coordinate
ring of an affine algebraic group. Then by the  Hochschild-Kostant-Rosenberg theorem
 $HH_\bullet(k \lbrack G \rbrack  ;k)\cong \land^{\bullet}(Lie(G))$ and hence is 
 independent of the group structure.

 Next we discuss the second point  above. Some interesting cyclic cocycles
were defined by Connes in the context of Lie algebra homology and group cohomology.
 For  example let $A$ be an algebra and $\delta _1,  \delta_2 :A \rightarrow A $ two 
 commuting derivations.
 Let $\tau  :A \rightarrow \mathbb{C}$  be an {\it{invariant trace }}  in the sense that $\tau$ 
  is a trace  and 
 $\tau(\delta_1(a))=\tau(\delta_2(a))=0$ for all $a \in A $. Then one can directly check 
 that the following is a cyclic $2$-cocycle on $A$~  \cite{ac80}
 :$$\varphi(a_0,a_1,a_2)=\tau(a_0(\delta_1(a_1)\delta_2(a_2)-\delta_2(a_1)\delta_1(a_2))).$$
 This cocycle is non-trivial. For example, if $A=A_\theta$ is the algebra of smooth noncommutative torus and
 $ e \in A_\theta$ is the smooth Rieffel projection, then $\varphi(e,e,e)=\pm q$, 
 where $\tau(e)=\mid p-q\theta\mid$~\cite{ac80}. 
 
  For a second example let $G$ be a (discrete) group and ${c}$ 
 be  a normalized group cocycle on $G$ with trivial coefficients. Then one can easily check that the following 
 is a cyclic cocycle on the group algebra $\mathbb{C} G$~\cite{achm00}
 \begin{align*}
 \varphi (g_0,g_1 \dots, g_n )=\left\{ \begin{matrix} {c}(g_1,g_2\dots,g_n)~~\text{if}~~g_0 g_1
  \dots g_n=1\\
 0 ~~~~~~~~~~~~\text{otherwise}\end{matrix}\right.
 \end{align*}  
   
 It is highly desirable to understand the origin of these formulas, put them in 
 a conceptual 
 context and generalize 
 them. For example we need to know in the case where a Lie algebra $\mathfrak{g}$ 
 acts by derivations on an algebra $A$, 
 $\mathfrak{g}\rightarrow Der (A)$,  if there is a map
 \[\gamma :H_\bullet (\mathfrak{g},\mathbb{C})\rightarrow HC^\bullet (A).\]  
	  
  Now let us indicate how the cohomology theory defined by Connes-Moscovici 
  ~\cite{achm98,achm99} and its dual 
  version in ~\cite{kr}  resolve  both issues.
  Let $\mathcal{H}$  be a Hopf algebra. Let $\delta$  be {\it character  } and
   $\sigma$  a 
 {\it group like } element on $\mathcal{H}$, i.e. $\delta:\mathcal{H}\rightarrow k$ 
  is  an
 algebra map and $\sigma:k \rightarrow \mathcal{H}  $ a coalgebra map. Following 
  \cite{achm98,achm99}, we say $(\delta,\sigma)$
is a {\it modular pair } if $\delta \sigma =id_k$ and a {\it  modular pair in involution } 
if,  in addition,  $(\sigma^{-1}\widetilde{S})^2=id_\mathcal{H} $ where  the {\it  
twisted antipode }  $\tilde{S}$ is de
fined by 
$$\widetilde{S}(h)=\sum_{(h)}\delta(h^{(1)}) S(h^{(2)}).$$
 Given $\mathcal{H}$, and $(\delta,\sigma )$, 
Connes-Moscovici define a cocyclic module  $\mathcal{H}_{(\delta,\sigma)}^\natural$
 as follows. 
Let   $\mathcal{H}_{(\delta,\sigma)}^{\natural,0}=k$ and 
$\mathcal{H}_{(\delta,\sigma)}^{\natural,n}=\mathcal{H}^{\otimes n}$ , $n\geq 1$.
The coface, codegeneracy and cyclic operators $\delta_i$, $\sigma_i$, $\tau$ are 
defined by 
\begin{eqnarray*}
\delta_0(h_1 \otimes \dots \otimes h_n)&=& 1_\mathcal{H} \otimes h_1 \otimes\dots
 \otimes h_n \\
\delta_i(h_1 \otimes\dots \otimes h_n)&=& h_1 \otimes \dots \otimes\Delta (h_i)
\otimes\dots \otimes h_n
  \;\;\text{for}\;\;1\leq i \leq n \\
\delta _{n+1}(h_1 \otimes\dots \otimes h_n)&=& h_1 \otimes\dots \otimes h_n 
\otimes \sigma \\
\sigma_i(h_1 \otimes\dots \otimes h_n)&=& h_1 \otimes\dots \otimes\epsilon
(h_{i+1})\otimes\dots \otimes h_n 
 \;\;\text{for}\;\;0 \leq i \leq n \\
\tau(h_1 \otimes\dots \otimes h_n )&=&\Delta^{n-1}\widetilde{S}(h_1)\cdot(h_2 \otimes\dots \otimes h_n \otimes \sigma). 
\end{eqnarray*} 	

These formulas were discovered in~\cite{achm98} and then  proved in full generality
 in \cite{achm99}. In ~\cite{cr}, M. Crainic 
gave an alternative approach based on Cuntz-Quillen formalism of cyclic homology ~\cite{cq}.    Note 
that  
the cosimplicial module $\mathcal{H}_{(\delta,\sigma)}^\natural$  
is the cosimplicial module associated to the coalgebra $\mathcal{H}$ with 
coefficients in $k$ via the 
unit map and $\sigma$.
The passage from the cyclic homology of (co)algebras  to the cyclic homology of Hopf 
algebras is remarkably similar to passage from 
de Rham cohomology to Lie  algebra cohomology. The key idea in both 
cases is
 {\it invariant cohomology}.

It is not difficult  to see that the above complex  is an  exact analogue
 of {\it invariant cohomology }
 in noncommutative geometry.
 In fact,  under the multiplication  map $\mathcal{H}\otimes \mathcal{H} \rightarrow \mathcal{H}$ 
 the coalgebra $\mathcal{H}$
  is an $\mathcal{H}$-module coalgebra. Let $\hat{ \mathcal{H}}_\natural$ be the cocyclic  module  of the 
  coalgebra $\mathcal{H}$. The cocyclic module $\hat{ \mathcal{H}}_\natural$  becomes 
  a  cocyclic $\mathcal{H}$-module 
   via the diagonal action $\mathcal{H}\otimes \hat{ \mathcal{H}}_\natural\rightarrow 
   \hat{ \mathcal{H}}_\natural $. We have   
   ${\hat{ \mathcal{H}_\natural}}^{\delta}=\mathcal{H}_{(\delta,1)}^\natural$  where 
   ${\hat{ \mathcal{H}_\natural}}^{\delta}$ is the space of $\delta$-coinvariants.

   The cohomology groups $HP^\bullet_{(\delta,\sigma)}(\mathcal{H})$  are  so far computed for the 
   following Hopf algebras. For 
   quantum  universal enveloping algebras no examples are known except for $U_q(sl_2)$  that we recall below.
   \begin{itemize}

\item [1.] If $\mathcal{H}=\mathcal{H}_n$ is the  Connes-Moscovici Hopf algebra, we have ~\cite{achm98} 
$$HP^n_{(\delta,1)}(\mathcal{H})\cong \bigoplus_{i=n\; (\text{mod}\;2)} H^i (\mathfrak{a}_n,\mathbb{C})$$
where $\mathfrak{a}_n$ is the Lie algebra of formal vector fields on $\mathbb{R}^n$.
\item[2.] If $\mathcal{H}=U(\mathfrak{g})$ is the enveloping algebra of a Lie algebra  $\mathfrak{g}$, we 
have ~\cite{achm98}
$$HP^n_{(\delta,1)}(\mathcal{H})\cong \bigoplus_{i=n\; (\text{mod}\; 2)} H_i(\mathfrak{g},\mathbb{C}_\delta)$$
\item[3.] If $\mathcal{H}=\mathbb{C}\lbrack G \rbrack $ is the coordinate ring of a nilpotent 
affine algebraic  group $G$, we have ~\cite{achm98}
$$HP^n_{(\epsilon,1)}(\mathcal{H})\cong \bigoplus_{i=n\;( \text{mod}\; 2)} H^i(\mathfrak{g},\mathbb{C}),$$  
where $\mathfrak{g}=Lie(G)$.
\item[4.] If $\mathcal{H}$ admits a normalized left Haar integral, then \cite{cr}
$$HP^1_{(\delta,\sigma)}(\mathcal{H})=0,\qquad HP^0_{(\delta,\sigma)}(\mathcal{H})=k.$$
Recall that a linear map $\int :\mathcal{H} \rightarrow k$ is called a normalized left  Haar integral 
if for all $h\in\mathcal{H}$,   $\int(h)=\int(h^{(1)})h^{(2)}$ and $\int(1)=1$. 
Compact quantum groups, finite dimensional Hopf algebras over a filed of characteristic zero,  
 and group  algebras are known to admit normalized Haar
 integral in the above sense. In the latter case $\int :kG \rightarrow k$
sending $g \mapsto 0$ for all $g\neq e$  and $e \mapsto 1$ is a Haar integral. Note that $G$ need not to be finite.
\item[5.] If $\mathcal{H}=U_q(sl_2(k))$ is the quantum universal algebra of $sl_2(k)$, we have ~\cite{cr},
$$HP^0_{(\epsilon,\sigma)}(\mathcal{H})=0,\quad  HP^1_{(\epsilon,\sigma)}(\mathcal{H})=k \oplus k.$$
\item[6.]Let $\mathcal{H}$ be a commutative Hopf algebra. The periodic cyclic cohomology of the cocyclic module
        $\mathcal{H}^\natural_{(\epsilon,1 )}$  can be  computed in terms of the Hochschild homology of coalgebra
         $\mathcal{H}$ with trivial 
        coefficients.

  \begin{prop}(\cite{kr})
  Let $\mathcal{H}$ be a commutative Hopf algebra. Its periodic cyclic 
  cohomology in the sense 
  of Connes-Moscovici is given by
  $$HP_{(\epsilon,1)}^n(\mathcal{H})=\bigoplus_{i=n\;(\text{mod}\;2)}H^i
  (\mathcal{H},k).$$
    \end{prop}
    For example, if $\mathcal{H}=k\lbrack G\rbrack$ is the algebra of 
    regular functions on
     an affine algebraic group $G$, the coalgebra complex of  
     $\mathcal{H}=k\lbrack G\rbrack$ is isomorphic 
     to the group cohomology complex of $G$ where  instead 
     of regular cochains one uses regular functions 
     $G\times G\times \dots \times G \rightarrow k$. Denote this cohomology by
      $H^i(G,k)$. It follows that 
   $$HP_{(\epsilon,1)}^n(k\lbrack G\rbrack)=\bigoplus_{i=n\;(\text{mod}\;2)}H^i(G,k).$$
   As is remarked in \cite{achm00}, if the Lie algebra Lie$(G)=\mathfrak{g}$ 
   is nilpotent, it  follows from 
   Van Est's theorem 
   that $H^i(G,k)\cong H^i(\mathfrak{g},k)$. This gives an alternative proof of
    Prop.4 and 
   Remark 5 in ~ \cite{achm00}.
   \end{itemize}

   Let $A$ be an $\mathcal{H}$-module algebra  and $Tr: A\rightarrow \mathbb{C}$  a $\delta$-{\it invariant 
   } linear map,  {{i.e.}}, $Tr(h(a))=\delta(h)Tr(a)$ for
 $h \in \mathcal{H},\;a \in A$. Equivalently,  $Tr$ satisfies the {\it integration by part } property:
 $$Tr(h(a)b)=Tr(a\tilde{S}(h)(b)).$$
 In addition we  assume $Tr(ab)=Tr(b\sigma a ).$
 Given $(A,\mathcal{H},Tr)$, Connes-Moscovici show that  the following map called the 
 {\it the characteristic map}, defines a morphism
 of cyclic modules $\gamma :\mathcal{H}^\natural_{\delta,\sigma}\rightarrow A^\natural$, 
 where $A^\natural=\hom(A_\natural,k)$ is the cocyclic module associated to $A$, 
  $$\gamma(h_1\otimes \dots\otimes h_n)(a_0,a_1,\dots,a_n)=Tr(a_0h_1(a_1)\dots h_n(h_n)). $$
 We therefore have  well-defined maps\\
 \begin{center}
 $\gamma : HC_{(\delta,\sigma)}^\bullet(\mathcal{H})\rightarrow HC^\bullet(A) $\\
  $\gamma : HP_{(\delta,\sigma)}^\bullet(\mathcal{H})\rightarrow HP^\bullet(A). $

 \end{center}
 
 Examples  show that, in general,   this map is non-trivial.  
   For example let $\mathfrak{g}$ be an abelian  $n$-dimensional Lie algebra 
   acting by derivations on an
   algebra $A$. Let $\delta_i\in Der(A)$ 
    be the  family of derivations corresponding to a basis $X_1$,\dots,$X_n$  
    of $\mathfrak{g}$, 
    and $Tr :A\rightarrow k$ an invariant 
    trace on $A$,  i.e. $Tr\delta_i(a)=0,\quad 1\le i\le n$. We have  
    $H_i(\mathfrak{g},k)\cong
     \land^i\mathfrak{g}$. In particular $H_n(\mathfrak{g},k)$ 
    is $1$-dimensional.  The inclusion $$H_n(\mathfrak{g},k)\hookrightarrow 
    \bigoplus_{i=n\;\text{mod}\;2}H_i(\mathfrak{g},k)\cong HP^n_{(\epsilon,1)}
    (U(\mathfrak{g}))$$
combined with the characteristic  map $\gamma$ defines  a map
 $$\gamma: H_n(\mathfrak{g},k)\cong k\rightarrow HC^n(A).$$
    The image of $X_1\land X_2\land \dots \land X_n$ under $\gamma $ is the cyclic
     $n$-cocycle $\varphi$ given by 
    $$\varphi(a_0,a_1,\dots,a_n)=\sum_{\sigma\in S_n}(-1)^nTr(a_0\delta_1
    (a_{\sigma(1)})\delta_2(a_{\sigma(2)})\dots 
    \delta_n(a_{\sigma(n)})).$$
   
   The rest of  this  section is devoted to a dual cyclic theory for Hopf 
   algebras which was defined, independently, 
   in \cite{kr,ta}.
   There is  a need for a dual theory to be developed. This is needed, for 
   example, when one studies
    coactions of Hopf algebras (or quantum 
   groups) on noncommutative spaces,
   since the original Connes-Moscovici
    theory works for actions only.
   A more serious problem is the fact that if $\mathcal{H}$ has normalized left  Haar integral 
   then its cyclic cohomology in the 
   sense of Connes-Moscovici    is trivial in positive dimensions ~\cite{cr},  but
    the dual theory is non-trivial. 
    
In \cite{kr} we  associated  a cyclic module to any Hopf algebra $\mathcal{H}$ 
over $k$ if $\mathcal{H}$ has 
a modular pair $( \delta , \sigma ) $ such that $\widehat{S}^2=id_\mathcal{H}$, where
 $\widehat{S}(h)=\delta(h^{(2)})\sigma S(h^{(1)})$. This cyclic module  can be 
 seen as the dual of the cocyclic module 
 introduced in $~\cite{achm99}$ 
by A. Connes and H. Moscovici. 
Using $\epsilon$ and $\delta $ one can endow $k$ with an $\mathcal{H}$-bimodule 
structure,{  i.e.},
$$\delta \otimes id : \mathcal{H} \otimes k \to k \quad and \quad  id \otimes 
\epsilon : k \otimes \mathcal{H} 
\to k. $$
Our cyclic module as a simplicial module is exactly the 
Hochschild complex of $ \mathcal{H} $ with coefficients in $ k $ where $ k $ is an 
 $\mathcal{H}$-bimodule as 
 above. So  if  we  denote our cyclic module by
${\widetilde {\mathcal{H}}}_\natural^{(\delta,\sigma)}, $  we have
 $ {\widetilde {\mathcal{H}}}_{\natural_n}^{( \delta,\sigma)}  = \mathcal{H}^
 {\otimes n} $, for $n> 0$  and $ {\widetilde {\mathcal{H}}}_{\natural_0}^{( \delta,\sigma)}
   = k $.   Its faces and degeneracies 
 are as follows: \\ 
\begin{eqnarray*}
 {\delta}_0 (h_1 \otimes h_2 \otimes  ... \otimes h_n) &=&\epsilon (h_1)h_2 \otimes h_3
  \otimes ... \otimes h_n \\
 {\delta}_i (h_1 \otimes h_2 \otimes 
 ... \otimes h_n) &=& h_1\otimes h_2 \otimes ...\otimes 
h_i h_{i+1}\otimes ... \otimes h_n   \\
{\delta}_n(h_1 \otimes h_2 \otimes  ... \otimes h_n) &=& \delta (h_n)h_1\otimes
 h_2 \otimes ... \otimes h_{n-1} \\
 {\sigma}_0(h_1 \otimes h_2 \otimes  ... \otimes h_n) &=&1 \otimes h_1  \otimes 
 ... \otimes h_n \\
\hspace{2cm}{\sigma}_i(h_1 \otimes h_2 \otimes ... \otimes h_n) &=&h_1 \otimes 
h_2 ...\otimes h_i \otimes 1 \otimes
 h_{i+1} ... \otimes h_n     \\
 {\sigma }_n(h_1 \otimes h_2 \otimes  ... \otimes h_n) &=& h_1 \otimes h_2  
 \otimes ...\otimes 1. 
\end{eqnarray*} 
To define a cyclic module it remains to introduce an action of cyclic 
group on our module. Our candidate is 
$$\tau _n(  h_1 \otimes h_2 \otimes  ... \otimes h_n ) = \sum  
 \delta(h_n^{(2)})\sigma S ( h_1^{(1)} h_2^{(1)} ...h_{n-1}^{(1)}h_n^{(1)} ) 
\otimes h_1^{(2)} \otimes ...\otimes h_{n-1}^{(2)} . $$

\begin {theorem}{(\cite{kr})}
Let $\mathcal{H} $ be a Hopf algebra over $k$  with a modular pair $(\delta , \sigma)$
 such that $\widehat{S}^2=id_\mathcal{H}$. 
Then $\widetilde {\mathcal{H}}_\natural^{( \delta,\sigma)}$
 with operators given   above   defines a cyclic module. Conversely, if $(\delta, 
 \sigma ) $ is a modular pair such that $\widetilde{\mathcal{H}}^{(\delta,\sigma)}_\natural$ is a cyclic module, then 
  $\widehat{S}^2=id_\mathcal{H}$. 
\end {theorem}       
Now let $A$ be an $\mathcal{H}$-comodule algebra. To define the characteristic 
map we need 
an analogue of an invariant trace.    
    
\begin{definition} A linear map, $Tr :A\rightarrow k$ is called $\delta$-trace  if 
\begin{equation*}
Tr(ab)= \sum _{(b)} Tr(b^{(0)}a)\delta (b^{(1)}) \hspace{2cm}\forall a , b\in  A.
\end{equation*} 
It is called $\sigma $-invariant if for  all $a , b\in  A$,
\begin{multline*}
 ~~~~~~~~~~~~~~~~~~~~~~~~~~~~~~~~~~~~~~~~~
 \sum _{(b)} Tr(ab^{(0)}) (b^{(1)})=\sum _{(a)}Tr(a^{(0)}b)S_{\sigma }(a^{(1)})\\ 
    \shoveleft{\text{or equivalently }}\\
    \shoveleft{~~~~~~~~~~~~~~~~~~~~ Tr(a^{(0)})a^{(1)}=Tr(a)\sigma.}\\
\end{multline*}
\end{definition}
Consider the map $\gamma :A_\natural \rightarrow \widetilde{\mathcal{H}}^
{(\delta,\sigma)}_\natural$ defined by 
$$\gamma(a_0\otimes a_1\otimes\dots \otimes a_n)=Tr(a_oa_1^{(0)}\dots a_n^{(0)})
a_1^{(1)}\otimes
  a_2^{(1)}\otimes\dots a_n^{(1)}.$$
It is proved in~\cite{kr} that $\gamma$ is a morphism of cyclic modules. 
    \begin{corollary}
Under the   above conditions, $ \gamma $ induces the following canonical maps: 

\begin{center}
$\gamma :\; {HC}_\bullet (A)\to \widetilde{HC}^{(\delta ,\sigma)}_\bullet
(\mathcal{H})$\\
$\gamma :\; {HP}_\bullet (A)\to \widetilde{HP}^{(\delta ,\sigma)}_\bullet
(\mathcal{H})$.
\end{center}
\end{corollary}

Next, we state a theorem which computes  the cyclic homology of cocommutative 
Hopf algebras. 
   \begin {theorem}(\cite{kr})\label{t1}
If $\mathcal{H}$ is a cocommutative Hopf algebra, then
$$ \widetilde{HC}^{({\delta,1})}_n(\mathcal{H}) = \bigoplus _{i\geq 0} H_{n-2i}
(\mathcal{H},   k_\delta),$$ where $k_\delta$ is the one dimensional module defined by $\delta$.

\end{theorem}
\begin{example}
Let $\mathfrak{g}$ be a Lie algebra over $k$ and $U(\mathfrak{g})$ be 
its enveloping algebra. One knows that 
$H_n(U(\mathfrak{g}) ; k) = H_n(\mathfrak{g} ; k)$ ~\cite{ld}.
 So by Theorem \textnormal{\ref{t1}} we have  
$$\widetilde{HC}^{(\delta,1)}_n(\mathfrak{g}) = \bigoplus _{i\geq 0}H_i(\mathfrak{g} ; k_\delta). $$
\end{example}
\begin{example}

Let $G$ be a discrete group and $\mathcal{H}=kG$ its group algebra. Then from theorem \textnormal{\ref{t1}} we have 
\begin{center}
$\widetilde{HC}_n^{(\epsilon,1)}(kG)\cong\bigoplus _{i\ge 0}H_{n-2i}(G,k)$\\
$\widetilde{HP}_n^{(\epsilon,1)}(kG)\cong\bigoplus _{i=n\;(\text{mod}\; 2)}H_{i}(G,k).$
\end{center} 
\end{example}
\begin{example}
Let $G$ be a discrete group and $\mathcal{H}=\mathbb{C}G$. Then the algebra 
$\mathcal{H}$ is a comodule algebra for 
the Hopf algebra $\mathcal{H}$ via coproduct map $\mathcal{H}\longrightarrow 
\mathcal{H}\otimes\mathcal{H}$. The map
 $Tr:\mathbb{C}G \rightarrow \mathbb{C}$  defined by $$Tr(g)=\left\{\begin{matrix}
  1 \quad\qquad g=e\\
0 \quad\qquad g \neq e \end{matrix}\right.$$  
is a $\delta$-invariant $\sigma$-trace for $\delta=\epsilon$, $\sigma=1$. 
The dual characteristic map\\ 
$\gamma^{\ast}:\widetilde{HC}^n_{(\epsilon,1)}(\mathbb{C}G)\rightarrow 
{HC}^n(\mathbb{C}G)$ 
 combined with the inclusion\\ 
$H^n(G,\mathbb{C})\hookrightarrow \widetilde{HC}^n_{(\epsilon,1)}
(\mathbb{C}G)$ is 
exactly the map $H^n(G,\mathbb{C})\rightarrow HC^n(\mathbb{C}G)$ 
described earlier in this section. 
\end{example}
It would be very interesting to compute the Hopf cyclic homology
$\widetilde{HC}_\bullet$ for compact 
quantum groups. Of course, one should look at algebra 
of polynomial or smooth functions on compact quantum groups, the 
$C^\ast$-completion
 being uninteresting from cyclic theory point of view. In the  following we recall two results
  that are known so far about
  quantum groups. 
  
   Let $k$ be a field of characteristic zero and $q\in k$, $q\neq 0$ and $q$ 
not a root of unity. The Hopf algebra $\mathcal{H}=A(SL_q(2,k))$ is defined as follows.
 As an algebra it is generated by symbols $a, b, c, d,$ with the following relations:
$$ ba=qab, \;\; ca=qac, \;\; db=qbd,\;\;dc=qcd,$$
$$bc=cb,\;\; ad-q^{-1}bc=da-qbc=1. $$
The coproduct, counit and antipode of   $\mathcal{H}$ are defined by  
$$\Delta (a)=a\otimes a+b\otimes c,\;\;\;\Delta (b)=a\otimes b+b\otimes c $$ 
$$\Delta (c)=c\otimes a+d\otimes c,\;\;\;\Delta (d)=c\otimes b+d\otimes d $$
$$\epsilon (a)=\epsilon (d)=1,\;\;\;\epsilon (b)=\epsilon (c)=0,$$
$$S(a)=d,\;\; S(d)=a,\;\;S(b)=-qb,\;\;S(c)=-q^{-1}c. $$
  For more details about $\mathcal{H}$ we refer to ~\cite{kl}. Because $S^2\neq id $, to define 
our cyclic structure we  need a modular pair $(\sigma, \delta)$ in involution. Let $\delta $ be as follows:
$$\delta (a)=q,\; \; \delta(b)=0,\; \; \delta (c)=0,\; \; \delta (d)=q^{-1}.$$
And $\sigma =1$. Then we have $\widetilde S_{(1, \delta )}^2=id$.\\

For computing cyclic homology we should at first compute the Hochschild homology
 $H_{\ast}(\mathcal{H}, k)$ where $k$ is an $\mathcal{H}$-bimodule via   $\delta$, $\epsilon$ for 
 left and right action of $\mathcal{H}$, respectively.\\
One knows $H_{\ast}(\mathcal{H}, k)=Tor_{\ast}^{\mathcal{H}^e}(\mathcal{H}, k)$, 
where $\mathcal{H}^e = \mathcal{H}\otimes \mathcal{H}^{op}$. So we need 
a resolution for $k$, or $\mathcal{H}$ as $\mathcal{H}^e$-module. We take advantage 
of the  free resolution for $\mathcal{H}$  given by Masuda et.al. \cite{kr}. By a 
lengthy computation one can check that 
$H_0(\mathcal{H},k)=0$, $H_1(\mathcal{H},k)=H_2(\mathcal{H},k)=k\oplus k$, and $H_n(\mathcal{H},k)=0$ for all 
$n\geq 3$. Moreover we find that  the operator
 $B=(1-\tau)\sigma N :H_1(\mathcal{H},k)\longrightarrow H_2(\mathcal{H},k)  $ is bijective and we obtain:

\begin{theorem}(\cite{kr})
For any $q\in k$ which is not a root  of unity one has\\ $\widetilde{ HC}_1(A(SL_q(2,k)))=k\oplus k$ and 
$\widetilde{HC}_n(A(SL_q(2,k)))=0$ for all $n\neq 1$. \\ In particular,  $\widetilde{HP}_0(A(SL_q(2,k)))
=\widetilde{HP}_1(A(SL_q(2,k)))=0$.  
\end{theorem}

The above theorem shows that Theorem $4.2$ is not true for non-cocommutative Hopf algebras.

The quantum universal enveloping algebra $U_q(sl(2,k))$ is an $k$-Hopf algebra
 which is generated as an $k$- algebra by symbols $\sigma$, $\sigma^{-1}$, $x$, $y$
subjecte to the following relations\\
$$\sigma \sigma^{-1}=\sigma^{-1} \sigma=1, \;\; \sigma x =q^2x\sigma,\;\;  \sigma y =q^{-2}y\sigma, \;\;
xy-yx = \frac{\sigma-\sigma^{-1}}{q-q^{-1}}.$$

The coproduct, counit and antipode of $U_q(sl(2,k))$ are defined by:
$$\Delta(x)=x\otimes\sigma +1\otimes x ,\;\; \Delta(y)=y\otimes 1+\sigma^{-1}\otimes y,
 \;\;\Delta(\sigma)=\sigma\otimes \sigma,$$
$$S(\sigma)=\sigma^{-1},\;\;S(x)=-x\sigma^{-1},\;\;S(y)=-\sigma
 y,$$$$\;\;\varepsilon(\sigma)=1,\;\varepsilon(x)=\varepsilon(y)=0.$$

It is easy to check thet $S^2(a)=\sigma a \sigma^{-1}$, so that  $(\sigma^{-1},
\varepsilon)$ is a modular pair in
 involution. As the first step to compute its cyclic homology we should find its 
Hochschild homology group with trivial coeficients. ($k$ is a $U_q(sl(2,k))$ 
bimodule via $\varepsilon$). We define a
 free resolution for $\mathcal{H}=U_q(sl(2,k))$ as a $\mathcal{H}^e$-module 
 as follows
\[
\begin{CD}
(*) \hspace{2cm} \mathcal{H}@<\mu<< M_0 @<d_0<< M_1 @< d_1<< M_2@< d_2<< M_3
 \dots
\end{CD}
\]
where  $M_0$ is $\mathcal{H}^e$, $M_1$ is  the free $\mathcal{H}^e$-module 
generated by symbols $1\otimes e_\sigma,
 1\otimes e_x, 1\otimes e_y$, $M_2$ is the free $\mathcal{H}^e$-module  
 generated by symbols $1\otimes e_x\land
 e_\sigma, 1\otimes e_y\land e_\sigma, 1\otimes e_x\land e_y $, and finally $M_3$ 
 is generated by $1\otimes e_x\land
   e_y\land e_\sigma$  as a free $\mathcal{H}^e$-module. We let $M_n=0$ for all $n\geq 4$. We claim that with the
    following boundary operators,   $(*)$ is a free resolution for $\mathcal{H}$
\begin{eqnarray*}
 &&d_0(1\otimes e_x)= x\otimes 1-1\otimes x \\
 &&d_0(1\otimes e_y)=  y\otimes 1-1\otimes y \\
 &&d_0(1\otimes e_\sigma)=  \sigma\otimes 1-1\otimes \sigma \\  
 &&d_1(1\otimes e_x\land e_\sigma)=(\sigma \otimes 1-1\otimes q^2\sigma)\otimes e_\sigma -(q^2x\otimes1-1\otimes x
  )\otimes e_x\\
&& d_1(1\otimes e_y\land e_\sigma)=(\sigma \otimes 1-1\otimes q^{-2}\sigma)\otimes e_\sigma -(q^{-2}y\otimes
 1-1\otimes y )\otimes e_y \\
&& d_1(1\otimes e_x\land e_y)=(y\otimes 1-1\otimes y)\otimes e_x-(x\otimes 1- 1\otimes x )\otimes e_y \\ 
  &&~~~~~~~~~~~~~~~~~~~  +\frac{1}{q-q^{-1}}(\sigma^{-1}\otimes \sigma^{-1}+ 1\otimes 1)\otimes e_\sigma \\
&&d_2(1\otimes e_x\land e_y\land e_\sigma )=(y\otimes 1-1\otimes q^2y)\otimes e_x\land e_\sigma\\
&&~~~~~~~~~ - q^2(q^2x\otimes 1-1\otimes x )\otimes e_y\land e_\sigma + q^2(\sigma\otimes 1-1\otimes \sigma)\otimes
 e_y\land e_x
\end{eqnarray*}

To show that this complex is a resolution, we need a homotopy map. First we recall that  the set
 \hbox{$\{\sigma^lx^my^n\mid l\in\mathbb{Z},m,n\in\mathbb{N}_0  \}$} is a P.B.W. type basis   for $\mathcal{H}$
  ~\cite{kl}.\\
Let $$\phi(a,b,n)=(a^{n-1}\otimes 1+a^{n-1}\otimes b\dots a\otimes b^{n-1}+1\otimes b^{n-1})$$  
 where $n\in \mathbb{N}, a\in \mathcal{H}, b\in \mathcal{H}^o $, and $\phi(a,b,0)=0$, and
  $\omega(p)=1 $ if $p \geq 0$
  and $0$ otherwise.\\
The following maps define a homotopy map  for $(*)$ i.e. $sd+ds=1$:
\begin{eqnarray*}
&&S_{-1}: \mathcal{H} \rightarrow M_0,\\
&&S(a)=1\otimes a, \\
&& S_0 : M_0 \rightarrow M_1,\\
&&S_0(\sigma^lx^my^n\otimes b)= (1\otimes b)((\sigma^lx^m\otimes 1)\phi(y,y,n)\otimes e_y+\\
&&~~~~~~~+(\sigma^l\otimes y^n)\phi(x,x,m)\otimes e_x)+\omega(l)(1\otimes x^my^n)
\phi(\sigma,\sigma,l)\otimes e_\sigma
 \\
&&~~~~~~~~~~~~~~~~~+(\omega(l)-1)(1\otimes x^my^n)\phi(\sigma^{-1},\sigma^{-1},-l)(\sigma^{-1}\otimes
 \sigma^{-1}\otimes e_\sigma),\\
&& S_1:M_1\rightarrow M_2,\\
&& S_1(\sigma^lx^my^n\otimes b \otimes e_y)=0,\\
&& S_1(\sigma^lx^my^n\otimes b \otimes e_x)=(1\otimes b)((\sigma^lx^m\otimes 1)\phi(y,y,n)\otimes e_x\land e_y\\ 
&&+\frac{1-q^{2n}}{(q-q^{-1})(1-q^2)}(\sigma^l\otimes y^{n-1})\phi(x,x,m)(\sigma^{-1}\otimes \sigma^{-1}+q^{-2}
\otimes 1 )\otimes e_x\land e_\sigma\\
  &&~~~~~~~~+\frac{1}{q-q^{-1}}(\sigma^lx^m\otimes 1)\phi(y,y,n-1)(\sigma^{-1}\otimes\sigma{-1}+q^2\otimes 1)\otimes
   e_y\land e_\sigma),\\
&& S_1(\sigma^lx^my^n\otimes b \otimes e_\sigma)=(1\otimes b)(q^2(\sigma^lx^m\otimes 1)\phi(y,q^2y,n)\otimes 
e_y\land e_\sigma\\
&&~~~~~~~~~~~~~~~~~~~~~~~~~~~~~~~~~~~~~~+q^{2(n-1)}(\sigma^l\otimes y^n)\phi(x,q^{-2}x,m)
\otimes e_x\land e_\sigma),\\
&& S_2:M_2\rightarrow M_3,\\
&&S_2(a\otimes b\otimes e_x\land e_y)=0,\\
&&S_2(a\otimes b \otimes e_y\land e_\sigma )=0,\\
&&S_2(\sigma^lx^my^n\otimes b\otimes e_x\land e_\sigma)=(1\otimes b)(\sigma^lx^m\otimes 1)
\phi(y,q^2y,n)\otimes
 e_x\land e_y\land e_\sigma,\\
&&S_n=0:M_n\rightarrow M_{n+1}~~\text{for}~~ n\geq 3.
\end{eqnarray*}

Again, by a rather long, but straight forward computation, we can check that $ds+sd=1$.    
By using the  definition of Hochschild homology as $Tor^{\mathcal{H}^e}(\mathcal{H},k)$ 
we have the following theorem:
\begin{theorem}(\cite{kr})
$H_0(U_q(sl(2,k)),k)=k$ and $H_n(U_q(sl(2,k)),k)=0$ for all $n\neq 0$ where $k$ is
 $U_q(sl(2,k))$-bimodule
 via $\varepsilon$ for 
both side. 
\end{theorem}
\begin{corollary}
$\widetilde{HC}_n(U_q(sl(2,k)))=k$ when $n$ is even, and  $0$ otherwise. 
\end{corollary}

 \section{ Cohomology of extended  Hopf algebras}
 In their study of  index theory for transversely elliptic operators and in order to treat the general
 non-flat case, Connes and Moscovici \cite{achm01} had to replace their Hopf algebra $\mathcal{H}_n$ 
 by a so-called \textquoteleft\textquoteleft extended Hopf algebra" $\mathcal{H}_{FM}$. In fact
  $\mathcal{H}_{FM}$ is neither a Hopf algebra 
 nor a Hopf algebroid in the 
sense of \cite{lu},
but it has enough structure to define a cocyclic module  similar to Hopf algebras 
~\cite{achm98, achm99, achm00}.\\

In attempting to define a cyclic cohomology theory for Hopf algebroids in general, 
we were led instead 
to define a closely related concept that we call an extended Hopf algebra. This 
terminology is already 
used in \cite{achm01}. All examples of  interest, including the Connes-Moscovici 
algebra $\mathcal{H}_{FM}$ 
are extended Hopf algebras.
    
Our first goal in this section is to recall the definition of {\it extended Hopf algebra} from \cite{kr2}
. This is closely related, but different from, {\it Hopf algebroids} in \cite{ lu,xu}. 
The reason we prefer this 
concept to Hopf
 algebroids is that it is not clear how to define cyclic homology of Hopf algebroids, 
 but it can be defined for 
 extended Hopf 
 algebras as we will recall from \cite{kr2}. The whole theory is motivated by \cite{achm01}.
 
 Broadly speaking, extended Hopf algebras and Hopf algebroids are quantizations
  (i.e. not necessarily commutative or cocommutative analogues) of 
  groupoids and Lie algebroids. This should be compared with the point of view that 
  Hopf algebras are quantization 
  of groups 
  and Lie algebras. Commutative Hopf algebroids were defined as cogroupoid objects 
  in the category of commutative algebras in \cite{ra}. The main example being  algebra of functions on a 
  groupoid. The 
  concept was later generalized to allow noncommutative total algebras. A decisive step was 
  taken in \cite{lu} where both 
  total and base algebra are allowed to be noncommutative. 
  
   To define a cocyclic  module  one  needs an 
   {\it antipode pair} $(S,\widetilde{S})$
    as we define below. Motivated by this observation and also the fundamental work of \cite{achm01},
   we were   led to define extended Hopf algebras and their cocyclic modules. 
   
      Recall from \cite{lu,xu} that a {\it bialgebroid } $(H,R,\Delta,\varepsilon)$ consist of 
    \begin{itemize}
\item[1:] An algebra $H$, an algebra $R$, an algebra 
homomorphism $\alpha :R \rightarrow H$, and  an algebra anti homomorphism $\beta : R \rightarrow H$ such that 
the image of $\alpha $ and $\beta $ commute in $H$. It follows that $H$ can be 
regarded as $R$-bimodule via
  $axb=\alpha(a)\beta(b)x$~~~~ $a,b\in R$~$x\in H.$\\
$H$ is called the ${total~ algebra}$, $R$ the ${base~ algebra }$, $\alpha$ the 
${source~ map}$ and $\beta$ 
the ${target ~map}$.   
\item[2:]A coproduct, i.e.  an $(R,R)$-bimodule map $\Delta : H\rightarrow H\otimes_RH$ 
with $\Delta(1)=1\otimes_R1$  satisfying the following conditions
\begin{itemize}
\item[i)]Coassociativity :
$$(\Delta \otimes_R id_H)\Delta = (id_H\otimes_R\Delta)\Delta :H \rightarrow H\otimes_RH\otimes_RH.$$
\item[ii)] Compatibility with product:
$$\Delta(a)(\beta(r)\otimes 1-1\otimes\alpha(r) )=0 \; \text{in} H\otimes_RH \;\text{for any}\;r\in R \;\;\;a\in H$$
$$\Delta (ab)=\Delta(a)\Delta(b)\;\;\text{for any }\;a,b\in H.$$
\end{itemize}
\item[3:] A counit, i.e. an $(R,R)$-bimodule  map $\epsilon:H\rightarrow R$ 
satisfying \\
$\epsilon(1_H)=1_R$ and $(\epsilon \otimes_Rid_H )\Delta =(id_H\otimes_R\epsilon)
\Delta =id_H:H\rightarrow  H.$ 
\end{itemize} 

\begin{definition}
Let $(H,R,\alpha , \beta, \Delta ,\varepsilon)$ be a $k$-bialgebroid. We call it
 a $Hopf$ $ algebroid$ if there 
is a bijective map $S:H\rightarrow H$ which is a antialgebra map satisfying the 
following conditions,
\begin{itemize}
\item[i)]\label{one}$S\beta = \alpha $.
\item[ii)]\label{three}$m_H(S\otimes id)\Delta =\beta \epsilon S:H\rightarrow H$.
\item[iii)]\label{three} There exists a linear map $\gamma:H\otimes_RH \rightarrow H\otimes H$ satisfying\\
$\pi\circ\gamma =id_{H\otimes_RH} :H\otimes_RH\rightarrow H\otimes_RH$ and $m_H(id\otimes S)\gamma 
\Delta =\alpha \epsilon :H\rightarrow H$\\
where $\pi:H \otimes H \rightarrow H \otimes _RH $ is the  natural projection. 
\end{itemize}  
\end{definition}

\begin{definition}\label{def}
Let $(H,R)$ be a bialgebroid. An {\rm antipode pair} $(S,\widetilde{S})$ consists of 
maps  $S,\widetilde{S}: H \rightarrow H$
 such that
 \begin{itemize}

\item [(i)] $S$ and $\widetilde{S}$ are antialgebra maps. 
\item[(ii)] $\widetilde{S}\beta=S \beta=\alpha.  $
\item[(iii)] $m_H(S \otimes id)\Delta =\beta \epsilon S:H \rightarrow H $ and 
 $m_H(\widetilde{S} \otimes id)\Delta =\beta \epsilon \widetilde{S}:H \rightarrow H.$
 
\item[(iv)] There exists a $k$-linear section $\gamma:  H\otimes_RH \longrightarrow H\otimes H$
 for the natural projection 
 $H\otimes H\longrightarrow H\otimes_RH$ such that  the map $\gamma\circ \Delta: H\longrightarrow H\ot H$ 
  is coassociative and the following two diagrams  are commutative: 
$$
 \xymatrix{H\ar[d]_\Delta\ar[r]^S & H\ar[r]^{\hspace{-5pt}\Delta} & H\otimes_R H \\
 H\otimes_RH\ar[d]_\gamma & &H\otimes H\ar[u]_\pi\\
 H\otimes H \ar[rr]^\tau & &  H\otimes  H \ar[u]_{S\otimes S} }
~~~~~~~~~\xymatrix{H\ar[d]_\Delta\ar[r]^{\widetilde{S} } & H\ar[r]^{\Delta} & H\otimes_R H \\
 H\otimes_RH\ar[d]_\gamma & &H\otimes H\ar[u]_\pi\\
 H\otimes H \ar[rr]^\tau & &  H\otimes  H \ar[u]_{S\otimes \widetilde{S} }}
$$
In the above diagrams  $\tau: H\otimes H\longrightarrow H\otimes H$ is the
 ``twisting map"  defined by $\tau(h_1\otimes h_2)=h_2\otimes h_1$.  Equivalently, 
 and by abusing the language, we say $S$ is an ``anticoalgebra map" and $\widetilde{S}$ is a 
 ``twisted anti coalgebra map", 
 i.e. for all $h\in H$ 
 \begin{equation}\label{anticoalgebra}
 \Delta S(h)=\sum S(\mathfrak{h}^{(2)})\otimes_RS(\mathfrak{h}^{(1)}),
 \end{equation}
 \begin{equation}\label{twisted anticoalgebra}
\Delta \widetilde{S}(h)=\sum S(\mathfrak{h}^{(2)})\otimes_R \tilde{S}(\mathfrak{h}^{(1)}),
\end{equation}
 where $\gamma(\Delta(h))=\sum\mathfrak{h}^{(1)}\otimes \mathfrak{h}^{(2)}.$
  \end{itemize}   
\end{definition}

\begin{definition}
An {\rm extended Hopf algebra} is a bialgebroid endowed with an antipode pair
 $(S,\widetilde{S})$ such that $\widetilde{S}^2=id_H$.
\end{definition}
\begin{bf}   
  Remark.
\end{bf}
The exchange operator $H \otimes_RH \rightarrow H \otimes _RH$, $x \otimes _R y \mapsto y \otimes _Rx$,
 is not well-defined in
 general.  A careful look at the proof of the cocyclic module property for the Connes-Moscovici cocyclic 
 module $\mathcal{H}^{(\delta,1)}_\natural$
  of a Hopf algebra $\mathcal{H}$ ( cf.  Theorem 2.1 in \cite{kr} )
 reveals that relations  (\ref{anticoalgebra}) and (\ref{twisted anticoalgebra}) (for $k=R$) play a fundamental role. 
 The same is true for Theorem  \ref{main}, but since $R$ is noncommutative in general, 
 these relations  make sense only after we fix a section $\gamma$  as  in Definition \ref{def}. 
 Coassociativity of the map $\gamma\circ\Delta: H\longrightarrow H\ot H$ is needed in the proof of Theorem 
 \ref{main}.  This motivates  our definition of an extended Hopf algebra. 
 
 Recall the Connes-Moscovici algebra ($\mathcal{H}_{FM},R$) 
 associated to a smooth manifold $M$ \cite{achm01}. It is shown  in \cite{achm01} that $\mathcal{H}_{FM}$ 
   is a free $R\otimes R$-module where $R=C^{\infty}(FM)$ is the algebra of smooth functions 
    on  the frame bundle $FM$. 
    In fact  fixing a torsion free connection on $FM$, one obtains a 
    Poincar\'e-Birkhoff-Witt type basis for $\mathcal{H}_{FM}$ over $R\otimes R$ 
    consisting of differential operators $Z_I\cdot\delta_K$, where 
    $Z_I$ is a  product of horizontal vector fields  $X_i$, $1\le i\le n$ and vertical vector fields $Y_j^i$ and $\delta_K$
    is a product of vector fields $\de$. 
    The coproduct $\Delta$ and the twisted antipode $\widetilde{S}$ are already 
    defined in \cite{achm01} and all the identities of a bialgebroid are verified. 
    All we have to do is to define a section  
    $\gamma: \mathcal{H}_{FM}\otimes_R \mathcal{H}_{FM}\longrightarrow \mathcal{H}_{FM}\otimes \mathcal{H}_{FM} $,
     an antipode  map $S: \mathcal{H}_{FM}\longrightarrow \mathcal{H}_{FM}$ and verify  
     the remaining  conditions of Definition \ref{def}. 
     
     To  this end, we first define $S$ on the generations of $\mathcal{H}_{FM}$ by 
  \begin{equation}\label{seh}
  \begin{matrix}
&S(\alpha(r))=\beta(r), && S(\beta(r))=\alpha(r),\\
 &S(Y_i^j)=-Y_i^j, &&S(X_k)=-X_k+\delta^i_{kj}Y^j_i,\\
  & S(\delta^i_{jk})=-\delta^i_{jk}.
  \end{matrix}
\end{equation}
We then extend $S$ as an antialgebra map, using the Poincar\'e-Birkhoff-Witt basis of $\mathcal{H}_{FM}.$ 

We define a section $\gamma: \mathcal{H}_{FM}\otimes_R\mathcal{H}_
{FM}\longrightarrow \mathcal{H}_{FM}\otimes \mathcal{H}_{FM}$ by the formula 
\begin{displaymath}
\gamma(\alpha(r)\otimes x\otimes \beta(s)\otimes_R\alpha(r')\otimes 
x'\otimes \beta(s'))=\alpha(r)\otimes x\otimes \beta(s)\alpha(r')\otimes 1\otimes x'\otimes \beta(s')
,\end{displaymath}
 where we use 
 the fact that $\mathcal{H}_{FM}$ is a free $R\otimes R$-module.   The following proposition is proved in \cite{kr2}.
\begin{prop}
The Connes-Moscovici algebra $\mathcal{H}_{FM}$ is an extended Hopf algebra.
\end{prop} 
   We give a few more examples of extended Hopf algebras.
\begin{itemize}
\item[1.]
Let $\mathcal{H}$ be a $k$-Hopf algebra, $\delta:\mathcal{H}\rightarrow k$ a character, i.e. 
 an algebra homomorphism  
and $\widetilde{S}_\delta=\delta\ast S$
 the $\delta$-twisted antipode defined by $\widetilde{S}_\delta(h)=\sum\delta(h^{(1)})S(h^{(2)})$, 
  as in \cite{achm98}. Assume that  $\widetilde{S}_\delta^2=id_{\mathcal{H}}$. 
 Then $(\mathcal{H},\alpha,
 \beta, \Delta,\epsilon,S,\widetilde{S}_\delta)$ is an  extended Hopf algebra, where $\alpha=\beta : k\longrightarrow 
 \mathcal{H} $  is  the unit map. More generally, given  any $k$-algebra $R$, let
  $H=R\otimes \mathcal{H}\otimes R^{op}$, where $R^{op}$ denotes the opposite algebra of $R$. 
 With the following structure $H$ is an  extended Hopf algebra over $R$: 
\begin{eqnarray*}
\alpha(a)&=&a\otimes 1\otimes1 \\
\beta(a)&=&1\otimes 1 \otimes a\\
\Delta(a\otimes h\otimes b)&=&\sum a\otimes h^{(1)}\otimes 1\otimes_R 1\otimes h^{(2)}\otimes b\\
\epsilon(a\otimes h\otimes b)&=&\epsilon(h)ab\\
S(a\otimes h\otimes b)&=&(b\otimes S(h)\otimes a) \\
\widetilde{S}(a\otimes h\otimes b)&=&(b\otimes \widetilde{S}_\delta(h)\otimes a), 
\end{eqnarray*}
 and the section $\gamma: H\otimes_R H \rightarrow H\otimes H$ is defined by 
 $\gamma(r\otimes h\otimes s\otimes_R r'\otimes h'\otimes s' )=r\otimes h\otimes sr' \otimes 1\otimes h'\otimes s'. $
  Then one can check that $(H,R)$ is an extended  Hopf algebra.
\item[2.]

The universal enveloping algebra $U(L,R)$  of a Lie-Rinehart algebra $(L,R)$ is an  extended Hopf algebra
 over the algebra $R$. For $X\in L$ and $r\in R$, we define\\   
$ 
\left. \begin{array}{cc}
\Delta({X})={X}\otimes_R 1+1\otimes_R{X} &\Delta({r})=  r\otimes_R 1 \\
\epsilon( X)=0 & \epsilon( r)=r\\
S( X)=- X & S( r)=  r.
\end{array}\right.$\\
Using  the Poincar\'e-Birkhoff-Witt theorem of Rinehart \cite{rin}, we extend $\Delta$ to be a multiplicative map, 
 $S$ to be  an antimultiplicative  map and $\epsilon$ by $\epsilon(rX_1\dots X_n)=0$ for $n\ge 1.$
 The source and target maps are the natural embeddings 
 $\alpha=\beta: R\hookrightarrow U(L,R)$ and $\widetilde{S}=S.$  
 We define the section $\gamma: U(L,R)\otimes_R U(L,R)\longrightarrow  U(L,R)\otimes U(L,R)$ by
 $ \gamma(rX_1\dots X_n\otimes_RsY_1\dots Y_m)=rsX_1\dots X_n\otimes Y_1\dots Y_m$. One can check that $\gamma $ is 
 well defined and $U(L,R)$ is an extended Hopf algebra. 
\item[3.]
Let $\mathcal{G}$ be a groupoid over a finite base (i.e., a category with a finite set of  objects, such that
 each morphism is invertible). Then the groupoid algebra $H=k\mathcal{G}$ is generated by 
morphism $g\in \mathcal{G}$ with unit $1=\sum_{X\in \mathcal{O}bj(\mathcal{G}) }id_X$, and the product  of two 
morphisms is equal to their composition  if the latter is defined and $0$ otherwise. It becomes an 
extended  Hopf algebra over $R=k\mathcal{S}$, where $\mathcal{S}$ is the subgroupoid of $\mathcal{G}$ 
whose objects  are those of
 $\mathcal{G}$ and $\mathcal{M}or(X,Y)=id_X$ whenever $X=Y$ and $\emptyset$ otherwise.
 The relevant  maps are defined for $g\in \mathcal{G}$ by $\alpha=\beta: R\hookrightarrow H$ is natural embedding, 
  $\Delta(g)=g\otimes_Rg $,   by $\epsilon(g)=id_{target(g)}$,   $S(g)=g^{-1}$.
The section  $\gamma : H\otimes_R H \longrightarrow H\otimes H  $ is defined by $\gamma(h\otimes_Rg)=h\otimes g$.
 It can easily be  checked that $H$ is both a Hopf algebroid and and extended Hopf algebra with $\widetilde{S}=S$
\end{itemize}

Given an extended Hopf algebra $(H,R)$ we define a cocyclic module $H_\natural$ as follows:
 $$H_\natural^0=R,\;\text{and}\;  H_\natural^n=H\otimes_R\otimes_R\dots\otimes_RH 
  \qquad(n\; \text{factors} ),\;\; n\geq 1.$$  
 The coface, codegeneracy and cyclic 
 actions $\delta_i$, $\sigma_i$ and $\tau$  are defined   by   
 \begin{eqnarray*}
\delta_0(h_1\otimes_R\dots \otimes_Rh_n)&=&1_H\otimes_R h_1\otimes_R\dots \otimes_Rh_n \\
\delta_i(h_1\otimes_R\dots \otimes_Rh_n)&=& h_1\otimes_R\dots\otimes_R\Delta (h_i)\otimes_R\dots \otimes_Rh_n
  \;\;\text{for}\;\;1\leq i\leq n \\
\delta_{n+1}(h_1\otimes_R\dots \otimes_Rh_n)&=&h_1\otimes_R\dots \otimes_Rh_n\otimes_R 1_H \\
\sigma_i(h_1\otimes_R\dots \otimes_Rh_n)&=& h_1\otimes_R\dots\otimes_R\epsilon(h_{i+1})\otimes_R\dots \otimes_Rh_n 
 \;\;\text{for}\;\;0\leq i\leq n \\
\tau(h_1\otimes_R\dots \otimes_Rh_n )&=&\Delta^{n-1}\widetilde{S}(h_1)\cdot(h_2\otimes\dots \otimes h_n\otimes 1_H). 
\end{eqnarray*}
These formulas were obtained in \cite{achm01} by transporting a cocyclic submodule of 
$A_\natural$ via a faithful trace to ${\mathcal{H}_{FM}}_\natural$, where $A$  is an algebra on 
which $\mathcal{H}_{FM}$ acts. In \cite{kr2}  we proved directly that these formulas  define a cocyclic  
modules for any extended Hopf algebra.
\begin{theorem}{\cite{kr2}}\label{main}
For any extended Hopf algebra $(H,R)$, the above formulas define a cocyclic module structure  on $H_\natural$. 
\end{theorem} 
The periodic cyclic cohomology of the universal enveloping algebra of Lie-Rinehart algebras is computed 
in ~\cite{kr2}.
Lie-Rinehart algebras interpolate between Lie algebras and commutative algebras, exactly in the same way that
 groupoids interpolate between groups and spaces. In fact Lie-Rinehart algebras can be considered as the 
 infinitesimal analogue of groupoids. In the following all algebras are unital, and all modules are unitary. 
  For more information on Lie-Rinehart algebras one can see  \cite{win,hub,rin}. 

Let $k$ be a commutative ring. A Lie-Rinehart algebra over $k$ is a pair $(L,R)$ where $R$ is a commutative 
$k$-algebra, $L$ is a $k$-Lie algebra and a left $R$- module, $L$ acts on $R$ by derivations 
$\rho: L \longrightarrow \mathcal{D}er_k(R) $ such that 
$\rho\lbrack X,Y\rbrack = \lbrack\rho(X),\rho(Y)\rbrack$ for all $X,Y$ in $L$ and the action is $R$-linear, 
and the Leibniz property holds:
$$\lbrack X,aY\rbrack =a\lbrack X,Y\rbrack +\rho(X)(a)Y\;\;\text{for all }\;X, Y \in L\; \text{and}\;a\in R. $$        
Instead of $\rho(X)(a)$ we simply write $X(a)$.
\begin{example}
Let $R=C^{\infty}(M)$ be the algebra of smooth functions on a manifold $M$ and 
$L=C^{\infty}(TM)=\mathcal{D}er_{\mathbb{R}}(C^{\infty}(M))$, the Lie algebra of vector fields on $M$.
Then $(L,R)$ is a Lie-Rinehart algebra, where the action $\rho:\; 
L=\mathcal{D}er_{\mathbb{R}}(R)\longrightarrow\mathcal{D}er_{\mathbb{R}}(R)  $ is the identity map.    
\end{example}
\begin{example}
Let $R=C^{\infty}(M)$ and $(L,R)$ a Lie-Rinehart algebra such that $L$ is a finitely generated 
projective $R$-module. Then it follows from  Swan's theorem that $L=C^{\infty}(E)$, is the space of
 smooth sections of a vector bundle over $M$. Since
  $\rho :C^{\infty}(E)\longrightarrow C^{\infty}(TM) $ is $R$-linear,  it is induced by a bundle 
  map $\rho : E\rightarrow TM. $\\
In this way we recover Lie algebroids as a particular example of Lie-Rinehart algebras.  
\end{example}
Next we recall the definition of the homology of a Lie-Rinehart algebra \cite{rin}.
This homology theory  is a simultaneous generalization of Lie algebra homology and 
de Rham homology. Let $(L,R)$ 
be a Lie-Rinehart algebra. A module over $(L,R)$ is a left $R$-module $M$ and a left Lie 
$L$-module $\varphi :L\rightarrow End_k(M)$,   denoted by $\varphi(X)(m)=X(m)$ such that 
for all $X\in L$, $a\in R$ and $m\in M$, 
\begin{center}
$X(am)=aX(m)+X(a)m$\\
$(aX)(m)=a(X(m)).$
\end{center}
  
 Alternatively, we can say an $(L,R)$-module is an $R$-module endowed with a $flat$ $connection$ 
  defined by $\nabla_X(m)=X(m)$,~ $X\in L,\;m\in M$.

Let $C_n=C_n(L,R;M)=M\otimes_R\mathcal{A}lt_R^n(L)$, where  $\mathcal{A}lt_R^n(L)$ 
denotes the $n$-th 
exterior power of the $R$-module $L$ over $R$. Let $d: C_n\longrightarrow C_{n-1}$ 
be the differential defined by \\
\begin{center}
$d(m\otimes X_1\land\dots \land X_n)=\sum_{i=1}^n(-1)^{i-1}X_i(m)\otimes X_1
\land\dots \land\hat{X_i}\dots\land X_n$\\
$+\sum_{1\leq i<j\leq n}(-1)^{i+j}m\otimes\lbrack X_i,X_j\rbrack\land X_1 \dots \land\hat{X_i}\dots
 \land\hat{X_j}\dots\land X_n.$
\end{center}
It is easy to check that $d^2=0$ and thus  we have  a complex $(C_n,d)$. The homology of this complex is,  
by definition, the homology of the Lie-Rinehart algebra $(L,R)$ with coefficients in $M$ and we denote this 
homology by $H_{\ast}(R,L;M)$.
To interpret this homology theory as a derived  functor, Rinehart
in \cite{rin} introduced the {\textit{universal enveloping
 algebra}} 
of a Lie-Rinehart algebra (L,R). It is an associative $k$-algebra, denoted ${U}(L,R)$, 
such that the category of $(L,R)$-modules as defined above is equivalent to the category of $U(L,R)$-modules. 
It is defined as follows.\\      
One can see easily that the following bracket defines a $k$-Lie algebra structure on $R\oplus L$: 
 $$\lbrack r+X,s+Y\rbrack =\lbrack X,Y\rbrack +X(s)-Y(r)\;\;\text{for}\; r,s\in R \;\text{and}\;X,Y \in L. $$ 
Let  $\tilde{U}=U(R\oplus L)$, be the  enveloping algebra of the Lie algebra $R\oplus L$, and let 
$\tilde{U}^+$ be the subalgebra generated by the canonical image of $R\oplus L$ in $U$. Then $U(L,R)=\tilde{U}^+/I
$, where $I$ is the two sided ideal generated by the set $\{(r.Z)'-r'Z' \mid  r\in R
 \;\text{and} \;Z\in R\oplus L\}$. 
 In ~\cite{rin}                      
Rinehart showed that if $L$ is a projective $R$-module, then \\
\begin{center}
$H_\ast(L,R;M)\cong Tor_\ast^{U(L,R)}(R,M).$\\
\end{center}     
  Next we compute the cyclic cohomology groups of the extended Hopf algebra $U(L,R)$ of a 
Lie-Rinehart algebra $(L,R)$. Let $S(L)$ be the symmetric algebra of the $R$-module $L$. 
It is an extended Hopf algebra over $R$. In fact  it is the enveloping algebra of the pair 
$(L,R)$ where $L$ is an abelian Lie algebra acting by zero derivations on $R$. Let $\land(L) $ 
be the exterior algebra of the $R$-module $L$. The following lemma computes the Hochschild cohomology 
of the cocyclic module $S(L)_{\natural}$.
\begin{lem}
Let $R$ be a commutative $k$-algebra and let $L$ be  a flat $R$-module. Then 
$$HH^{\ast}(S(L)_{\natural})\cong \land^{\ast}(L).$$
\end{lem}
The following proposition computes the periodic cyclic cohomology of the extended Hopf algebra  $U(L,R)$ 
associated to a Lie-Rinehart algebra $(L,R)$ in terms of its Rinehart homology. It extends a similar result 
for enveloping algebra of Lie algebras from \cite{achm98}.
\begin{prop}{(\cite{kr2})}
If $L$ is a projective $R$-module,  then we have 
$$HP^n(U(L,R))=\bigoplus_{i= n \;\text{mod} \;2}H_i(L,R;R),$$
where $HP^*$ means periodic cyclic cohomology.
\end{prop}
\begin{corollary}
Let $M$ be a smooth closed manifold and $\mathcal{D}$ be the algebra of differential operators on $M$.
It is an extended Hopf algebra and its periodic cyclic homology is given by 
$$HP_n(\mathcal{D})=\bigoplus_{i=n \;\;(\text{mod}\; 2)}H^i_{dR}(M).$$
\end{corollary}
\begin{proof}
We have $\mathcal{D}=U(L,R)$, where $L=C^{\infty}(TM)$ and $R=C^\infty(M)$. Dualizing the above proposition,
 we obtain 
$$HP_n(\mathcal{D})=\bigoplus_{i=n\;(\text{mod}\;2)}H^i(L,R)=\bigoplus_{i=n\;(\text{mod}\;2)}H^i_
{dR}(M)$$ 
\end{proof}
\begin{definition}{(Haar system for bialgebroids )} Let $(H,R)$ be a bialgebroid.\\
Let $\tau : H\longrightarrow R$ be a  right R-module map. We call $\tau$ a left Haar system for $H$ if 
$$\sum_{(h)}\alpha(\tau(h^{(1)})) h^{(2)}=\alpha(\tau(h))1_H$$
and $\alpha \tau=\beta\tau$.
We call $\tau$ a normal left Haar system if $\tau(1_H)=1_R$.
\end{definition}
We give a few examples of Haar systems. Let $H$  be the Hopf algebroid of a groupoid with finite 
base. Then it is easy to see that $\tau :H\rightarrow R $ defined by $\tau(id_x)=id_x$  for
all $x\in Obj(\mathcal{G})$ and $0$ otherwise is a normal Haar system for $H$. 
 In a related example, one can directly  check that  the map $\tau : A_\theta \rightarrow \mathbb{C}\lbrack  
 U,U^{-1}\rbrack$
defined by $$\tau(U^nV^m)=\delta_{m,0}U^n$$ is a normal Haar system for the noncommutative torus $A_\theta $.
 It is shwn in \cite{kr2}  that  $A_\theta$ is an extended Hopf algebra over $\mathbb{C}\lbrack U,U^{-1}\rbrack$. 

\begin{prop}
Let  $H$ be an extended Hopf algebra that admits a normal left Haar system. Then 
$HC^{2i+1}(H)=0$ and $HC^{2i}(H)=\ker\{\alpha-\beta\}$ for all $i\geq 0$.
\end{prop}
Finally in this section we compute  the Hopf periodic cyclic cohomology of commutative Hopf 
algebroids in terms of Hochschild cohomology.  
Given an extended Hopf algebra $(H,R)$, we denote the Hochschild cohomology of the cocyclic module 
$H_\natural$ by $H^i(H,R)$.
 It is the cohomology of the complex \\
$$
\begin{CD}
R@>d_0>> H@>d_1>> H\otimes_RH @>d_2>> H\otimes_RH\otimes_RH@>d_3>>\dots  
\end{CD}
$$
Where  the first differential  is $d_0=\alpha -\beta $ and $d_n $ is given by 
\begin{multline*}
d_n(h_1 \otimes_R \dots \otimes_R h_n)=1_H \otimes_R h_1 \otimes_R \dots \otimes_R h_n +\\
\sum_{i=1}^n(-1)^i h_1 \otimes_R \dots\otimes_R \Delta (h_i)\otimes_R \dots \otimes_R h_n  + \\
(-1)^{n+1}  h_1 \otimes_R \dots \otimes_R h_n \otimes_R 1_H. 
\end{multline*}

By a commutative extended Hopf algebra we mean an extended Hopf algebra $(H,R)$ where  both $H$ and $R$  are 
 commutative 
 algebras. In \cite{kr}, it is shown that  the periodic cyclic cohomology, in the sense of Connes-Moscovici, of a 
 commutative Hopf algebra admits a simple description. In fact, if $\mathcal{H}$ is a commutative Hopf algebra 
 then we have (\cite{kr}, Theorem 4.2):
 \begin{equation}\label{kr2}
 HC^n_{(\epsilon,1)}(\mathcal{H})\cong \bigoplus_{i\ge 0}H^{n-2i}(\mathcal{H},k),
 \end{equation}
 where the cohomologies on the right hand side are Hochschild cohomology of the coalgebra $\mathcal{H}$ with trivial 
 coefficients. Since the cocyclic  module of Theorem \ref{main}
  reduces to Connes-Moscovici cocyclic module if $H$ happens to be a Hopf algebra, 
 it is natural to expect that the analogue of isomorphism (\ref{kr2}) hold true for
  commutative extended  Hopf algebras. Furthermore, the analogue of Propositions 4.2 
  and 4.3 in \cite{kr},   which are crucial in  establishing 
   the above isomorphism (\ref{kr2}), are true for   extended Hopf algebras with similar proofs \cite{kr2}. This 
    leads us to the following conjecture. 

\begin{Conjecture}\label{conj} 
Let $(H,R)$ be a commutative extended Hopf algebra.
Then its  cyclic cohomology is given by
$$HC^n(H)\cong \bigoplus_{i\ge 0} H^{n-2i}(H,R).$$
$~~~\shoveright{\square}$ 
\end{Conjecture}
\section{Cohomology of smash products}
 A celebrated problem in cyclic homology theory is to compute the cyclic homology of the crossed product 
 algebra $A\ltimes G $, where the group $G$ acts on the algebra $A$ by automorphisms. If $G$   
 is a discrete group, there is a spectral sequence, due to Feigin and Tsygan~ \cite{ft}, which converges to the
  cyclic homology 
 of the crossed product algebra. This result generalizes
  Burghelea's calculation of the cyclic homology of a group algebra \cite{ld}.
  In \cite{gj}  Getzler and Jones gave a new proof of this spectral sequence using their Eilenberg-Zilber theorem 
  for cylindrical modules.
  In \cite{ak}, this spectral sequence has been extended to all Hopf algebras with invertible antipode.
   In this section we recall this result.
  
  Let $\mathcal{H}$ be a Hopf algebra and $A$ an $\mathcal{H}$-module algebra. We define a bicomplex, in fact a
   cylindrical module  $A\natural \mathcal{H}$ as follows: Let  
    $$(A\natural H){p,q}=\mathcal{H}^{\otimes(p+1)}\otimes A^{\otimes(q+1)} \qquad p,q\geq 0. $$   The vertical and 
    horizontal operators,    $\tau^{p,q}$, $\delta^{p,q}$, $\sigma^{p,q}$ and $t^{p,q}$,$d^{p,q}$,  $s^{p,q}$
     are defined by    
\begin{multline*}
{\tau^{p,q}( g_0 , \dots , g_p \mid a_0 , \dots , a_q) }
= (g_0^{(1)}, \dots , g_p^{(1)} \mid S^{-1} ( g_0^{(0)} g_1^{(0)} \dots g_p^{(0)}) \cdot a_q , a_0, \dots ,a_{q-1}) \\
\shoveleft{
\delta^{p,q}_i(g_0,\dots,g_p \mid a_0 , \dots , a_q)= (g_0,\dots,g_p \mid a_0,\dots , a_i a_{i+1},\dots , a_q) 
\;\;\; 0 \le i <q } \\
\shoveleft{\delta^{p,q}_q (g_0,\dots,g_p \mid a_0 , \dots , a_q)= (g_0^{(1)},\dots,g_p^{(1)} \mid (S^{-1}(g_0^{(0)}
 g_1^{(0)} \dots g_p^{(0)} ) 
\cdot a_q)a_0,\dots , a_{q-1})} \\
\shoveleft{\sigma^{p,q}_i(g_0,\dots,g_p \mid a_0 , \dots , a_q) = (g_0,\dots,g_p \mid a_0,\dots , a_i , 1 ,
 a_{i+1},\dots , a_q) 
\;\;\; 0 \le i \le q }\\
\shoveleft{t^{p,q}( g_0 , \dots , g_p \mid a_0 , \dots , a_q)
=(g_p^{(q+1)},g_0, \dots , g_{p-1} \mid g_p^{(0)} \cdot a_0, \dots
,g_p^{(q)} \cdot a_{q})} \\
\shoveleft{d^{p,q}_i(g_0,\dots,g_p \mid a_0 , \dots , a_q)
=(g_0,\dots,g_i g_{i+1},\dots,g_p \mid a_0,\dots , a_q) \;\;\; 0 \le i <q } \\ 
\shoveleft{d^{p,q}_q (g_0,\dots,g_p \mid a_0 , \dots , a_q)
= (g_p^{(q+1)}g_0,g_1,\dots,g_{p-1} \mid g_p^{(0)} 
\cdot a_0,\dots , g_p^{(q)} \cdot a_{q})}\\
\shoveleft{s^{p,q}_i (g_0,\dots,g_p \mid a_0 , \dots , a_q)
=(g_0,\dots,g_i,1,g_{i+1},\dots,g_p \mid a_0,\dots , a_q) \;\;\; 0 \le i \le q. } 
\end{multline*} 
\begin{bf}
Remark.
\end{bf} The cylindrical module $A\natural \mathcal{H}$ in ~\cite{ak} is defined for all Hopf algebras.
 For applications,  however, one  has to assume that $S$ is invertible. The above formulas are essentially 
 isomorphic to those in ~\cite{ak}, when $S$ is invertible.

\begin{theorem}{(\cite{ak})}
Endowed with the above operations, $A\natural \mathcal{H}$ is a cylindrical module.\\
$~~~\shoveright{\square } $
\end{theorem}
\begin{corollary}
The diagonal $d(A \natural \mathcal{H})$ is a cyclic module. \\
$~~~\shoveright{\square } $
\end{corollary}
Our next task is to identify the diagonal $d(A \natural \mathcal{H})$ with the cyclic module of the smash product 
$(A\# \mathcal{H})_\natural$. Define a map
 $\phi :(A\# \mathcal{H})_\natural\rightarrow d(A \natural \mathcal{H})$ by 
 \begin{eqnarray*}
 &&\phi ( a_0 \otimes g_0,\dots ,a_n \otimes g_n )=\\
& & (g_0^{(1)},g_1^{(2)},\dots,g_n^{(n+1)} \mid  S^{-1}(g_0^{(0)}g_1^{(1)} \dots g_n^{(n)}) 
\cdot a_0,S^{-1}(g_1^{(0)}g_2^{(1)}
 \dots g_n^{(n-1)}) \cdot a_1,\dots \\
& & \hspace{5cm} ,
S^{-1}(g_{n-1}^{(0)}g_n^{(1)})\cdot a_{n-1},S^{-1}(g_n^{(0)}) \cdot a_n)   
\end{eqnarray*}
By a long computation one shows that $\phi$ is a morphism of cyclic modules \cite{ak}.
\begin{theorem}{(\cite{ak})}
We have an isomorphism of cyclic modules  $d(A \natural \mathcal{H})\cong (A\# \mathcal{H})_\natural$.
\end{theorem}
\begin{proof}
Define a map $\psi :d(A\natural \mathcal{H})\rightarrow (A \# \mathcal{H})_\natural$ by 
$\psi (g_0, \dots , g_n \mid a_0 , \dots ,a_n)=$  
\begin{eqnarray*}
((g_0^{(0)} g_1^{(0)} \dots g_n^{(0)}) \cdot a_0 \otimes g_0^{(1)} ,(g_1^{(1)} \dots g_n^{(1)}) \cdot a_1
 \otimes g_1^{(2)},\dots , g_n^{(n)} \cdot a_n \otimes g_n^{(n+1)}).
\end{eqnarray*}
Then one can check that $\phi \circ \psi = \psi \circ \phi =id$.
\end{proof}   

Now we are ready to give an spectral sequence  to compute the cyclic homology of the smash product $A\#\mathcal{H}$.
 By using the Eilenberg-Zilber theorem for cylindrical modules, we have:
 \begin{theorem} \label{th:Apei-zi}
There is a
quasi-isomorphism of mixed complexes \[
\text{Tot} (  (A \natural \mathcal{H})) \cong 
 d (A \natural \mathcal{H} )\cong (A\#\mathcal{H})^\natural,\] 
and therefore an isomorphism of  cyclic homology groups,  \[ HC_{\bullet}(
\text{Tot} ( A \natural \mathcal{H})) \cong HC _{\bullet} ( A  \# H ). \]  
$~~~~\shoveright{~~~\square}$
\end{theorem}   
 Next, we show that one can identify the  $E^2$-term of the spectral sequence obtained 
 from the column 
 filtration. To this end, we define an action  of 
 $\mathcal{H}$  on the first row of $A \natural \mathcal{H}$, denoted by 
 $A_\mathcal{H}^\natural=\{ \mathcal{H}\otimes A^{\otimes(n+1)}\}_{n\geq 0} $
  by 
  \begin{equation*} 
h \cdot (g \mid a_0, \dots , a_n ) = (h^{(n+1)} \cdot g \mid h^{(0)} \cdot a_0 , \dots , h^{(n)} \cdot a_n) 
\end{equation*}
where $h^{(n+1)} \cdot g = h^{(n+1)} g \; S^{-1}(h^{(n+2)}) $ is an action of $\mathcal{H}$ on itself.
 We let $C_{\bullet}^{\mathcal{H}}(A)$ be  the space of  coinvariants  of $\mathcal{H} \otimes A^{\otimes (n+1)}$ 
under the above action. So in $C_{\bullet}^{\mathcal{H}}(A)$, we have 
\[h \cdot (g \mid a_0, \dots , a_n ) = \epsilon(h) (g \mid a_0, \dots , a_n ). \]
We define the following operators on $C_{\bullet}^{\mathcal{H}}(A)$, 
\begin{multline*}
~~~~~~~~~~~~~~~~~~~~~\tau_n(g \mid a_0 , \dots , a_n )= ( g^{(1)} \mid (S^{-1}
(g^{(0)}) \cdot a_n),a_0,\dots,a_{n-1})  \\
\shoveleft{~~~~~~~~~~~~~~~~~~~~~\delta_i(g \mid a_0 , \dots , a_n )= ( g \mid a_0,
\dots ,a_i a_{i+1},\dots,a_{n})} \\
\shoveleft{~~~~~~~~~~~~~~~~~~~~~\delta_n(g \mid a_0 , \dots , a_n )= ( g^{(1)} \mid (S^{-1}(g^{(0)})
 \cdot a_n)a_0,a_1,\dots,a_{n-1})}\\
\shoveleft{~~~~~~~~~~~~~~~~~~~~~\sigma_i(g \mid a_0 , \dots , a_n )= ( g \mid a_0,
\dots ,a_i,1,a_{i+1},\dots,a_{n})} \\
\end{multline*}
\begin{prop}{(\cite{ak})}
$C_{\bullet}^{\mathcal{H}}(A)$ with the operators defined above  is a cyclic module.
$~~\shoveright{\square}$
\end{prop}
Let $M$ be a left $\mathcal{H}$-module. Then $M$ is an $H$-bimodule if we let 
$\mathcal{H}$ act on the right on  $M$ via 
the counit map: $m.h=\varepsilon(h)m$. We denote the resulting  Hochschild homology 
groups by $H_\bullet(\mathcal{H},M)$. 
Explicitly    it is computed from the complex $C_p(\mathcal{H},M)=\mathcal{H}^{\otimes p}\otimes M, \quad p\geq 0$, 
with the differential $\delta :C_p(\mathcal{H},M)\rightarrow C_ {p-1}(\mathcal{H},M)$  defined by 
\begin{eqnarray*} 
\delta (g_1,g_2,\dots , g_p,m)= 
 \epsilon (g_1) (g_2,\dots ,g_p,m) +
\sum_{i=1}^{p-1} (-1)^i (g_1 , \dots ,g_i g_{i+1},\dots , g_p ,m) + \\(-1)^p (g_1,\dots,g_{p-1},g_p \cdot m). 
\end{eqnarray*}
Let $C_q(A_\mathcal{H}^\natural)=\mathcal{H}^{\otimes q}\otimes A_\mathcal{H}^\natural$ and let $\mathcal{H}$
 act on it by 
$h \cdot (g_1, \dots , g_p \mid m) = (g_1, \dots , g_p \mid h \cdot m)$, where the action of $\mathcal{H}$  on 
$A_\mathcal{H}^\natural$ is given by conjugation. So we can construct $H_p(\mathcal{H},C_q(A_\mathcal{H}^\natural))$.

Now we can show that our original cylindrical complex $(A \natural\mathcal{H},(\delta,\sigma,\tau),(d,s,t))$ 
can be identified with the cylindrical 
complex $(\mathsf{C}_p(\mathcal{H},\mathsf{C}_q ( A^{\natural}_{\mathcal{H}}),
( \mathfrak{d},\mathfrak{s},\mathfrak{t}),(\bar{\mathfrak{d}},\bar{\mathfrak{s}},
\bar{\mathfrak{t}}))$ under 
the transformations  $\beta : (A \natural \mathcal{H})_{p,q} \rightarrow 
\mathsf{C}_p(\mathcal{H},\mathsf{C}_q(A_{\mathcal{H}}^{\natural}))$
 and $\gamma :  \mathsf{C}_p(\mathcal{H},\mathsf{C}_q(A_{\mathcal{H}}^{\natural})) 
 \rightarrow  (A \natural \mathcal{H})_{p,q}$defined by 
 \begin{eqnarray*}
& &\beta ( g_0, \dots , g_p \mid a_0 , \dots , a_q ) = ( g_1^{(0)}, \dots , g_p^{(0)}
 \mid g_0 g_1^{(1)} \dots g_p^{(1)} \mid
a_0, \dots , a_q ) \\
& &\gamma ( g_1 , \dots , g_p \mid g \mid a_0 , \dots , a_q) = (g S^{-1}(g_1^{(1)} 
\dots g_p^{(1)}),g_1^{(0)}, \dots , g_p^{(0)}  \mid 
a_0, \dots , a_q ).
\end{eqnarray*}
One checks that $\beta\gamma =\gamma\beta=id.$
 To compute  the homologies of the mixed complex 
 $(Tot({C}(A \natural \mathcal{H}), b + \bar{b} + \mathbf{u} (B + \bar{B}))$
 we filter it by the subcomplexes (column filteration) 
 \[
\mathsf{F}^i_{pq} = \sum_{q \le i} (\mathcal{H}^{\otimes(p+1)} \otimes A^{\otimes(q+1)}) . \]

\begin{theorem}{(\cite{ak})}
The $\mathsf{E}^0$-term of the spectral sequence is isomorphic to the complex \[ 
\mathsf{E}^0_{pq}=(\mathsf{C}_p(\mathcal{H},\mathsf{C}_q(A^{\natural}_{\mathcal{H}})
  ),\delta) \]
and the $\mathsf{E}^1$-term is  \[
\mathsf{E}^1_{pq}= (H_p( \mathcal{H} , \mathsf{C}_q(A^{\natural}_{\mathcal{H}})  ) ,
 \mathfrak{b} + \mathbf{u} \mathfrak{B})). \]
The $\mathsf{E}^2$-term of the spectral sequence is
\[ \mathsf{E}^2_{pq} = HC_q( H_p(\mathcal{H}, \mathsf{C}_q(A^{\natural}_{\mathcal{H}}))), \]
the cyclic homologies of the cyclic module $H_p(\mathcal{H},\mathsf{C}_q(A^{\natural}_{\mathcal{H}}))$.
\end{theorem}

\section{Invariant Cyclic Homology}\label{sec7}    
 
 In this section,   we first  define the concept of  { \it Hopf triple}
  and its {\it invariant cyclic homology }.   One can think of invariant cyclic homology as  noncommutative 
  analogue of invariant  de Rham cohomology  as defined by Chevalley and Eilenberg  \cite{ce}.
   We indicate that cyclic homology of Hopf algebras is an example  of invariant  cyclic homology.
    We also present  our Morita invariance  theorem for invariant cyclic homology. Note that the
     result could not be formulated 
    for cyclic homology of Hopf algebras since the algebra of $n\times n$ matrices over 
    a Hopf algebra is not a Hopf algebra.  
    One can find the details of this section in ~\cite{kr5}. 
     In the second part of this section   we define the  invariant cyclic cohomology 
     of Hopf cotriples. One example is the 
Connes-Moscovici cyclic cohomology  of a Hopf algebra with a modular pair in involution 
in the sense of \cite{achm99} which turns  out to be the 
invariant cyclic cohomology  of the coalgebra $\mathcal{H}$. 
This  is implicit in \cite{achm98} and explicitly done in \cite{cr} for $\sigma=1$. We go, however, beyond this 
( fundamental) example and define  a cocyclic module for  any Hopf cotriple $(C,\mathcal{H},V)$ consisting  of  
an $\mathcal{H}$-module coalgebra $C$, an $\mathcal{H}$-comodule $V$ and a compatible character $\delta$ on
 $\mathcal{H}$.  Again the details can be found in \cite{kr5}.     
 \begin{definition}
 By a left  {\it Hopf triple } we mean a  triple $(A, \mathcal{H}, M)$, where $\mathcal{H}$ is a Hopf algebra,
  $A$ is a left $\mathcal{H}$\nobreakdash-comodule algebra and $M$ is a left $\mathcal{H}$\nobreakdash-module.
  Right Hopf triples are defined in a similar way.     
\end{definition}  
\begin{example}\label{example3.2}~~~
\begin{itemize}
\item [(i)](Trivial triples). Let $\mathcal{H}=k$, $M$ any $k$\nobreakdash-module, and $A$ any $k$-algebra. 
Then $(A,k,M)$ is a left Hopf triple.
\item[(ii)] Let $\mathcal{H}$ be a Hopf algebra and $M$ a left $\mathcal{H}$\nobreakdash-module.  Then 
$(A,\mathcal{H},M)$ is a left Hopf triple,  where $A=\mathcal{H}$ is the underlying algebra of $\mathcal{H}$
 and $\mathcal{H}$ coacts on $\mathcal{H}$ via its comultiplication. In particular, for $M=k$ and $\mathcal{H}$     
 acting on $k$ via a character $\delta$, we obtain a Hopf triple $(\mathcal{H},\mathcal{H},k_\delta)$.
 \end{itemize}
 \end{example}
  
   Given a left Hopf triple $(A, \mathcal{H},M)$, let $C_n(A,M)=M\ot A^{\ot(n+1)}$. 
   We define  simplicial and cyclic operators  on $\{C_n(A,M)\}_n$ by 
  \begin{align}\label{hasan}
 &\delta_0(m\ot a_0\ot a_1\ot\dots\ot a_n)=m\ot a_0a_1\ot a_2\ot\dots \ot a_n,\notag\\
&\delta_i(m\ot a_0\ot a_1\ot\dots\ot a_n)=m\ot a_0\ot \dots \ot a_i a_{i+1}\ot\dots \ot  a_n, &&1\le i\le n-1,\notag\\
& \delta_n(m\ot a_0\ot a_1\ot\dots\ot a_n)=a_n^{(-1)}m\ot a_n^{(0)}a_0\ot a_1\dots \ot a_{n-1},\\        
 &\sigma_i(m\ot a_0\ot a_1\ot\dots\ot a_n)=m\ot a_0 \ot \dots \ot a_i\ot 1 \ot \dots \ot a_n, &&0\le i\le n,\notag\\  
 &\tau(m\ot a_0\ot a_1\ot\dots\ot a_n)=a_n^{(-1)}m\ot a_n^{(0)}\ot a_0\ot\dots \ot a_{n-1}.\notag
 \end{align}
 One can check that endowed with the above operators,  $\{C_n(A,M)\}_{n}$ is a paracyclic module. 
 
 Next, we define a left  $\mathcal{H}$-coaction $\rho:C_n(A,M)\longrightarrow \mathcal{H}\ot C_n(A,M)$ by 
 $$\rho(m\ot a_0\ot\dots \ot a_n)=(a_0^{(-1)}\dots a_n^{(-1)})\ot  m\ot a_0^{(0)}\ot \dots \ot a_n^{(0)}.$$
With the above coaction,  $C_n(A,M)$ is an   $\mathcal{H}$\nobreakdash-comodule.
 To define the space of coinvariants  of $C_n(A,M)$, we fix a grouplike element $\sigma\in\mathcal{H}$.
Let
$$C_n^\mathcal{H}(A,M)= C_n(A,M)^{\text{co}\mathcal{H}}=\{ x\in C_n(A,M)\mid \rho(x)=\sigma\ot x\},$$
be the space of coinvariants  of $C_n(A,M)$ with respect to $\sigma$. We would like to find conditions that 
guarantee $\{C_n^\mathcal{H}(A,M)\}_n$ is a cyclic module. This leads us to the following definitions and results.
 
 \begin{definition}\label{definition3.5}
 Let $M$ be a left $\mathcal{H}$\nobreakdash-module and $\sigma\in \mathcal{H}$ a grouplike element.
  We define the $(M,\sigma)$\nobreakdash-{\it twisted antipode} $\widehat{S}:M\ot 
  \mathcal{H}\longrightarrow M\ot\mathcal{H}$ by
  $$\widehat{S}(m\ot h)=h^{(2)}m\ot \sigma S(h^{(1)}).$$
 \end{definition}
 \begin{definition}\label{def3.4}
 Let $M$ be a left $\mathcal{H}$\nobreakdash-module and $\sigma\in\mathcal{H}$ a grouplike element. We call
  $(M,\sigma)$ a 
 matched pair if $\sigma m=m$ for all $m\in M$. We call the matched pair $(M,\sigma)$ 
 a matched pair in involution if 
  $$(\widehat{S})^2=id: M\ot\mathcal{H}\longrightarrow M\ot \mathcal{H},$$
  where $\widehat{S}$ is defined in Definition \textnormal{{\ref{definition3.5}}}.
 \end{definition}
 \begin{example} \label{example3.6}
Let $M=k_\delta$ be  the one dimensional module defined by a character $\delta\in \mathcal{H}$. It is 
 clear that $(M,\sigma)$ is a matched pair in involution if and only if $(\delta, \sigma)$
  is a {\it modular pair in involution } 
  in the sense of  \cite{kr}, i.e.,  $\delta(\sigma)=1$ and $(\sigma\widetilde{S}_\delta)^2=id$.
  \end{example}
  \begin{lem}\label{lem3.10}
 Let $(A, \mathcal{H},  M)$ be a $\sigma$-compatible left Hopf triple. Then for any $a\in A$ and $m\in M$ 
 $$a^{(-1)}\sigma S(a^{(-3)})\ot a^{(-2)}m\ot a^{(0)}=\sigma\ot a^{(-1)}m\ot a^{(0)}.$$
 \end{lem}
  The following theorem is the main result that enable us to define the
   invariant cyclic homology of Hopf triples. 
  \begin{theorem}(\cite{kr5})\label{1}
 Let $(A, \mathcal{H}, M)$ be a  Hopf triple such that  $(M,\sigma)$ is a matched  pair in involution. 
 Then $\{ C_n^\mathcal{H}(A,M)\}_{n}$  endowed with  simplicial and 
 cyclic operators induced by \text{(1)}, is a 
 cyclic module. 
 \end{theorem}
 \begin{proof}
 As a first step we show that the induced simplicial and cyclic  operators are well
  defined on $\{C_n^\mathcal{H}(A,M)\}_n$ . We just  prove this for  $\tau$, and  $\delta_n$ 
  and leave the rest to the reader. Let $(m\ot a_0\ot\dots \ot a_n)\in C_n^\mathcal{H}(A,M)$. We have 
  \begin{equation}
     a_0^{(-1)}\dots a_n^{(-1)}\ot m\ot a_0^{(0)}\ot \dots a_n^{(0)}=\sigma\ot m\ot a_0\ot \dots a_n 
     \end{equation}
  which implies 
 $$
 a_0^{(-1)}\dots a_{n-1}^{(-1)}\ot m\ot a_0^{(0)}\ot \dots  \ot a_{n-1}^{(0)}\ot a_n=
  \sigma S(a_n^{(-1)})\ot m\ot a_0\ot \dots a_n^{(0)}$$
  and 
  \begin{multline*}
    a_n^{(-1)}a_0^{(-1)}\dots a_{n-1}^{(-1)}\ot a_n^{(-2)}
  m\ot a_n^{(0)}\ot a_0^{(0)}\ot \dots \ot a_{n-1}^{(0)}=\\
  \shoveright{~~~~=\underbrace{a_n^{(-1)}\sigma S(a_n^{(-3)})\ot a_n^{(-2)}m\ot 
  a_n^{(0)}}\ot a_0\ot \dots \ot a_{n-1}.}\\
  \end{multline*}
  Applying   Lemma \ref{lem3.10} for $a=a_n$ we have 
  \begin{multline}
    a_n^{(-1)}a_0^{(-1)}\dots a_{n-1}^{(-1)}\ot a_n^{(-2)}
  m\ot a_n^{(0)}\ot a_0^{(0)}\ot \dots \ot a_{n-1}^{(0)}=\\
  \sigma \ot a_n^{(-1)}m\ot a_n^{(0)}\ot a_0\ot \dots \ot a_{n-1}
  \end{multline}
  which means  $\tau(m\ot a_0\ot\dots \ot a_n)\in C_n^\mathcal{H}(A,M)$.\\ 
   
    From ($3$) we obtain  
    \begin{multline*}
    a_n^{(-1)}a_0^{(-1)}\dots a_{n-1}^{(-1)}\ot a_n^{(-2)}
  m\ot a_n^{(0)} a_0^{(0)}\ot \dots \ot a_{n-1}^{(0)}=\\
  \sigma \ot a_n^{(-1)}m\ot a_n^{(0)} a_0\ot \dots \ot a_{n-1}
  \end{multline*}
   which implies  $d_n(m\ot a_0\ot\dots \ot a_n)\in C_{n-1}^\mathcal{H}(A,M)$.\\
   Checking  that the  other simplicial operators  are well defined on $\{C_n^\mathcal{H}(A,M)\}_n$  is
    straightforward. 
    
   The only thing left is to show that $\tau^{n+1}=id$. We have    
   $$\tau^{n+1}(m\ot a_0\ot\dots \ot a_n)=a_0^{(-1)}\dots a_n^{(-1)} m\ot a_0^{(0)}\ot \dots a_n^{(0)}.$$   
   Now since we are in $C_n^\mathcal{H}(A,M)$,  and by ($2$) we have  
   $$\tau^{n+1}(m\ot a_0\ot\dots \ot a_n)=\sigma m\ot a_0\ot\dots a_n=m\ot a_0\ot\dots \ot a_n,$$
    because $(\sigma, M)$ is a matched pair.
   \end{proof}
   
   We denote the resulting Hochschild, cyclic and periodic cyclic homology groups 
   of the cyclic module $\{C_n^\mathcal{H}(A,M)\}_n$ by $HH_\bullet^\mathcal{H}(A,M)$,
    $HC_\bullet^\mathcal{H}(A,M)$ and
   $HP^\mathcal{H}_\bullet(A,M)$, respectively, and refer to them  as {\it invariant
    Hochschild}, {\it cyclic} and 
   {\it periodic cyclic } homology groups of the $\sigma$-compatible Hopf triple $(A,\mathcal{H},M)$.
   
   We give a few  examples of invariant cyclic homology. 
   More examples can be  found in \cite{kr5}.  
   It is clear that if $(A,k,k)$ is a trivial Hopf triple 
   ( Example \ref{example3.2} (i)), then $HC^k_\bullet(A,k)\cong HC_\bullet(A)$,
    i.e.,  in this case,  invariant cyclic homology is  the same  as cyclic homology of algebras.

	   We  show  that  cyclic homology of Hopf algebras in the sense of \cite{kr}  is an 
   example of invariant cyclic homology theory defined in this section. Consider 
   the $\sigma$-compatible Hopf triple 
   $(\mathcal{H},\mathcal{H},k_\delta)$ defined in Example \ref{example3.2}(ii). One can check  that for $M=k_\delta$ 
   the operators in \textnormal{(\ref{hasan})} are exactly the operators defined in \cite{kr}. This proves the following
    proposition. 
   \begin{prop}\label{prop3.14}
   The cyclic modules $\{\widetilde{\mathcal{H}}_n^{(\delta,\sigma)}\}_{n}$
    and $\{C_n^\mathcal{H}(\mathcal{H},k_\delta)\}_n$ are isomorphic. 
   \end{prop} 
 Let $(A, \mathcal{H},M)$ be a    Hopf triple.  
 One can easily see that  $(M_n({A}),\mathcal{H},M)$ is also a  Hopf triple, where
 the coaction of $\mathcal{H}$ on $M_n({A})$ is induced by the coaction  of $A$, i.e., 
 for all $a\ot u\in A\ot M_n(k)=M_n(A)$,  
{$\rho(a\ot u)=a^{(-1)}\ot a^{(0)}\ot u$}. We have the following theorem: 
  \begin{theorem}{(Morita invariance,\cite{kr5})}
  For any matched pair in involution $(M,\sigma)$,   and any  $k\ge 1$ one has 
  $$HC^\mathcal{H}_n(A,M)\cong HC_n^\mathcal{H}(M_k(A,M),\qquad n\ge 0.$$
    \end{theorem}
 Now let us  go to the dual case: coalgebras and {\it Hopf cotriples}. In this part as  promised 
  we define the cyclic cohomology of Hopf cotriples which  is a unification of  cyclic cohomology
   of coalgebras and cohomology of Hopf algebras introduced by  Connes-Moscovici in \cite{achm98}.      
\begin{definition}
By a left { Hopf cotriple} we mean a triple $(C,\mathcal{H},V)$ where $\mathcal{H}$ is a 
Hopf algebra, $C$ is a left $\mathcal{H}$\nobreakdash-module coalgebra and $V$ is 
a left  $\mathcal{H}$\nobreakdash-comodule.
\end{definition}
\begin{example}
Let $\mathcal{H}$ be a Hopf algebra and $V$ 
 a left  $\mathcal{H}$\nobreakdash-comodule. Then $(C,\mathcal{H},V)$ is a left Hopf cotriple,
 where $C=\mathcal{H}$ is the underlying coalgebra of $\mathcal{H}$, with  $\mathcal{H}$ acting on 
 $C$ via multiplication. In particular for $V=k_\sigma$ and $\mathcal{H}$ coacting on $k=k_\sigma$  via a grouplike 
 element $\sigma\in \mathcal{H}$, we have a Hopf cotriple $(\mathcal{H},\mathcal{H},k_\sigma)$. This is 
 the Hopf cotriple that is relevant to Connes-Moscovici theory \cite{achm98,achm99,achm00}. 
\end{example}
\begin{example}{(Trivial cotriples).} Let $C$ be a 
coalgebra, $\mathcal{H}=k$ and $V=k$. Then $(C,k,k)$ is a Hopf cotriple.
\end{example}

Given a Hopf cotriple $(C,\mathcal{H},V)$, let $C^n(C,V)=V\otimes C^{\ot(n+1)}$.
 We define cosimplicial and cyclic operators on $\{ C^n(C,V)\}_n=\{V\otimes C^{\ot(n+1)} \}_n$ by
\begin{eqnarray*}\label{f1}
\delta_i(v\otimes c_0 \otimes c_1 \otimes \dots \otimes c_n)&=& v\otimes c_0 \otimes\dots 
\otimes  c_i^{(1)}\otimes c_i^{(2)}\otimes c_n ~~~ 0 \leq i \leq n\\
\delta_{n+1}(v\otimes c_0 \otimes c_1 \otimes \dots \otimes c_n)&=& v^{({0})}\otimes c_0^{(2)}\otimes c_1 
\otimes \dots \otimes c_n  \otimes v^{({-1})}  c_0^{(1)}\\
\sigma_i(v\otimes c_0 \otimes c_1 \otimes \dots \otimes c_n)&=&v\otimes c_0 \otimes \dots 
c_i \otimes \varepsilon(c_{i+1})\otimes 
\dots\otimes c_n ~~~0\leq i \leq n-1\\
\tau(v\otimes c_0 \otimes c_1 \otimes \dots \otimes c_n)  &=&v^{({0})}\otimes c_1 \otimes c_2 \otimes 
\dots \otimes c_n \otimes v^{({-1})}c_0.
\end{eqnarray*}
One can check that endowed with the above  operators, $\{C^n(C,V)\}_{n}$ is a paracocyclic module.
 We have a diagonal $\mathcal{H}$-action on $C^n(C,V)$, defined by 
  $$h(v\ot c_o\ot\dots \ot c_n)=v\ot h^{(1)}c_0\ot h^{(2)}c_1\ot \dots \ot h^{(n+1)}c_n. $$
   It is easy to see that $C^n(C,V)$ is  an $\mathcal{H}$\nobreakdash-module.
   
   To define the space of coinvariants, we fix a character  of $\mathcal{H}$, say  $\delta$. Let
   $$C_\mathcal{H}^n(C,V)=\frac{ C^n(C,V)}{\text{span}
   \{hm-\delta(h)m\mid m\in C^n(C,V),\; h\in\mathcal{H}\}}$$
 be the space of coinvariants of $C^n(C,V)$ under the action of $\mathcal{H}$ and with respect
 to $\delta$.
 Our first task is to find sufficient  conditions under which  $\{C_\mathcal{H}^n(C,V)\}_n$ 
 is a cocyclic module. 
 
   Let us recall the twisted antipode $\widetilde{S}:\mathcal{H}\rightarrow \mathcal{H}$, where 
  $\widetilde{S}(h)=\delta(h^{(1)})S(h^{(2)})$, from ~\cite{achm98}.
  We define  the $V$-{\it twisted antipode}
  \begin{center}
  $\widetilde{S}_V:V\otimes\mathcal{H}\longrightarrow V\otimes\mathcal{H}$\\
  \shoveleft {by }\\
   $\widetilde{S}_V(v\otimes h)=v^{({0})}\otimes S^{-1}(v^{({-1})})\widetilde{S}(h). $
    \end{center}

    \begin{definition}
    We call the pair $(\delta,V)$ a { comatched pair} if 
    $$v^{({0})}\delta(v^{({-1})})=v \;\;\;\text{for all}\; v\in V.$$ 
     We call the comatched pair $(\delta, V)$ a { comatched pair in involution} if
     $$(\widetilde{S}_V)^2=id_{V\otimes\mathcal{H}}.$$ 
    \end{definition}
    
        Now we can state some of the main results of invariant cyclic 
        cohomology of Hopf cotriples. For more details see \cite{kr5}. 
    \begin{theorem}{\label{theorem4.10}}
    Let $(C,\mathcal{H},V)$ be a  Hopf cotriple such that $(\delta, V)$ is a comatched pair in involution. 
    Then $\{C^n_\mathcal{H}(C,V)\}_n$
     is a cocyclic module. 
   \end{theorem}
    \begin{example}
    Let $\mathcal{H}=V=k$ Then  for any coalgebra $C$ one has $\{C_k^n(C,k)\}_n$ is 
    the natural cyclic module, $C_\natural$, of the coalgebra $C$.
    \end{example}
    \begin{example}\label{lemma4.12}
    Let $\mathcal{H}$ be a Hopf algebra,  $C=\mathcal{H}$, and $V=k_\sigma$. The Hopf cotriple 
    $(\mathcal{H},\mathcal{H},k_\sigma)$ is $\delta$-compatible if and only if $(\delta,\sigma)$
     is a modular pair in involution 
    in the sense of ~\cite{achm99}. In this case $\{ C^n_\mathcal{H}(C,V)\}_n$ is isomorphic 
    to the Connes-Moscovici cocyclic module $\mathcal{H}^{(\delta,\sigma)}_\natural$.
    \end{example}

\end{document}